\documentclass[12pt]{article}
\usepackage{amsmath,amsthm,amssymb,graphicx,color}
\DeclareMathRadical{\sqrtsign}{symbols}{"70}{largesymbols}{"70}
\newcommand{\bb}{\mathbb}
\newcommand{\gothic}{\mathfrak}
\DeclareMathOperator{\Shadow}{Sh} \DeclareMathOperator{\Sh}{Sh}
\newcommand{\half}{{\bb H}}

\newcommand{\natls}{{\bb N}}
\newcommand{\ratls}{{\bb Q}}
\newcommand{\reals}{{\bb R}}

\newlength{\figboxwidth}
\newcommand{\makefig}[3]{
        \begin{figure}[htb]
        \refstepcounter{figure}
        \label{#2}
        \begin{center}
                #3~\\
                \smallskip
                Figure \thefigure.  #1
        \end{center}
        \medskip
        \end{figure}
}

\renewcommand{\bold}[1]{\medskip \noindent {\bf #1 }\nopagebreak}

\newcommand{\inv}{^{-1}}
\newcommand{\cross}{\times}
\newcommand{\st}{\;\: : \;\:}         %Such that

\newcommand{\Vol}{\operatorname{Vol}}

\newcommand\G{\Gamma}

\newcommand\Ra{\mathbb R}

\newcommand\Za{\mathbb Z}
\newcommand\Qa{\mathbb Q}
\newcommand\Na{\mathbb N}

\DeclareMathOperator{\SOL}{Sol} \DeclareMathOperator{\DL}{DL}
\def\@ifundefined#1#2#3%
  {\expandafter\ifx\csname#1\endcsname\relax#2\else#3\fi}

\@ifundefined{theoremstyle}{
}{
\theoremstyle{plain} %default
}
\newtheorem{theorem}{Theorem}[section]

\newtheorem{proposition}[theorem]{Proposition}
\newtheorem{lemma}[theorem]{Lemma}

\newtheorem{corollary}[theorem]{Corollary}

\@ifundefined{theoremstyle}{
}{
\theoremstyle{definition} %default
}
\newtheorem{definition}[theorem]{Definition}

\newcommand{\cC}{{\cal C}}
\newcommand{\cD}{{\cal D}}
\newcommand{\cE}{{\cal E}}
\newcommand{\cF}{{\cal F}}
\newcommand{\cG}{{\cal G}}

\newcommand{\cU}{{\cal U}}
\newcommand{\cV}{{\cal V}}

\mathchardef\GG="321D
\newcommand{\gb}{{\gothic b}}

\theoremstyle{plain}

\newcommand{\Sol}{\operatorname{Sol}}

\newcommand{\Area}{\operatorname{Area}}
\newcommand{\pSol}{\operatorname{Sol}(n')}

\begin{document}

\title{Coarse differentiation of quasi-isometries II: Rigidity for Sol and Lamplighter groups}
\author{Alex Eskin, David Fisher and Kevin Whyte}
%\thanks{First author partially supported by NSF grant DMS-0244542.
%Second author partially supported by NSF grants DMS-0226121 and
%DMS-0541917. Third author partially supported by NSF grant
%DMS-0349290 and a Sloan Foundation Fellowship.} }
\date{$ $}

\maketitle

\section{Introduction and statements of results}
\label{section:rigidity}

This paper continues the work announced in \cite{EFW0} and begun
in \cite{EFW1}.  For a more detailed introduction, we refer the
reader to those papers. As discussed in those papers, all our
theorems stated above are proved using a new technique, which we
call {\em coarse differentiation}. Even though quasi-isometries
have no local structure and conventional derivatives do not make
sense, we essentially construct a ``coarse derivative" that models
the large scale behavior of the quasi-isometry.  From this point
of view, the coarse derivatives of maps studied here are
constructed in \cite{EFW1} and this paper consists entirely of a
coarse analysis of coarsely differentiable maps.

We now state the main results whose proofs are begun in
\cite{EFW1} and finished here. The group $\SOL\cong
{\Ra}{\ltimes}\Ra^2$ with $\Ra$ acting on $\Ra^2$ via the diagonal
matrix with entries $e^{z/2}$ and $e^{-z/2}$.  As matrices, $\SOL$
can be written as :
\begin{displaymath}
\Sol = \left\{\left. \begin{pmatrix} e^{z/2} & x & 0 \\ 0 & 1 & 0 \\
0 & y & e^{-z/2}
    \end{pmatrix} \right| (x,y,z) \in \Ra^3 \right\}
\end{displaymath}
\noindent The metric $e^{-z}dx^2+e^{z}dy^2 + dz^2$ is a left
invariant metric on $\SOL$.  Any group of the form
${\Za}{\ltimes}_T\Za^2$ for $T \in SL(2,\Za)$ with $|tr(T)|>2$ is
a cocompact lattice in $\SOL$.

The following theorem proves a conjecture of Farb and Mosher:

\begin{theorem}
\label{theorem:sol} Let $\G$ be a finitely generated group
quasi-isometric to $\SOL$.  Then $\G$ is virtually a lattice in
$\SOL$.
\end{theorem}

We also prove rigidity results for wreath products $\Za{\wr}F$
where $F$ is a finite group. The name lamplighter comes from the
description $\Za{\wr}F = F^\Za \rtimes \Za$ where the $\Za$ action
is by a shift. The subgroup $F^\Za$ is thought of as the states of
a line of lamps, each of which has $|F|$ states. The "lamplighter"
moves along this line of lamps (the $\Za$ action) and can change
the state of the lamp at her current position. The Cayley graphs
for the generating sets $F \cup \{\pm 1\}$ depend only on $|F|$,
not the structure of $F$. Furthermore, $\Za {\wr} F_1$ and $\Za
{\wr} F_2$ are quasi-isometric whenever there is a $d$ so that
$|F_1| = d^s$ and $|F_2|=d^t$ for some $s,t$ in $\Za$.    The
problem of classifying these groups up to quasi-isometry, and in
particular, the question of whether the $2$ and $3$ state
lamplighter groups are quasi-isometric, were well known open
problems in the field, see \cite{dlH}.

\begin{theorem}
\label{theorem:lamplighter2} The lamplighter groups $\Za{\wr}F$
and $\Za{\wr}F'$ are quasi-isometric if and only if there exist
positive integers $d,s,r$ such that $|F|=d^s$ and  $|F'|=d^r$.
\end{theorem}
\noindent For a rigidity theorem for lamplighter groups, see
Theorem~\ref{theorem:lamplighter} below.

To state Theorem \ref{theorem:lamplighter} we need to describe a
class of graphs. These are the Diestel-Leader graphs, $DL(m,n)$,
which can be defined as follows: let $T_1$ and $T_2$ be regular
trees of valence $m+1$ and $n+1$. Choose orientations on the edges
of $T_1$ and $T_2$ so each vertex has $n$ (resp. $m$) edges
pointing away from it. This is equivalent to choosing ends on
these trees. We can view these orientations at defining height
functions $f_1$ and $f_2$ on the trees (the Busemann functions for
the chosen ends).   If one places the point at infinity
determining $f_1$ at the top of the page and the point at infinity
determining $f_2$ at the bottom of the page, then the trees can be
drawn as:

\begin{figure}[ht]
\begin{center}
\includegraphics[width=5.5in]{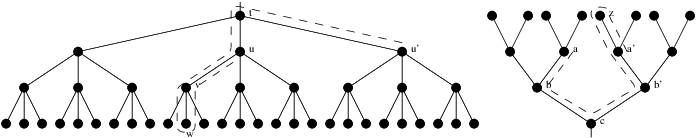}
\caption{The trees for $DL(3,2)$. Figure borrowed from
\cite{PPS}.} \label{fig:1}
\end{center}
\end{figure}

\noindent The graph $DL(m,n)$ is the subset of the product $T_1
\times T_2$ defined by $f_1 + f_2 = 0$. The analogy with the
geometry of $\SOL$ is made clear in \cite[Section 3]{EFW1}. For
$n=m$ the Diestel-Leader graphs arise as Cayley graphs of
lamplighter groups $\Za {\wr} F$ for $|F|=n$.  This observation
was apparently first made by R.Moeller and P.Neumann \cite{MN} and
is described explicitly, from two slightly different points of
view, in \cite{Wo2} and \cite{W}. We prove the following:

\begin{theorem}
\label{theorem:lamplighter} Let $\G$ be a finitely generated group
quasi-isometric to the lamplighter group $\Za{\wr}F$.  Then there
exists positive integers $d,s,r$ such that $d^s=|F|^r$ and an
isometric, proper, cocompact action of a finite index subgroup of
$\G$ on the Diestel-Leader graph $DL(d,d)$.
\end{theorem}

\noindent{\bf Remark:} The theorem can be reinterpreted as saying
that any group quasi-isometric to $DL(|F|,|F|)$ is virtually a
cocompact lattice in the isometry group of $DL(d,d)$ where $d$ is
as above.
\medskip

{\small {\it Acknowledgements.} First author partially supported
by NSF grant DMS-0244542. Second author partially supported by NSF
grants DMS-0226121 and DMS-0541917. Third author partially
supported by NSF grant DMS-0349290 and a Sloan Foundation
Fellowship. The second author thanks the University of Chicago for
hospitality while this paper was being written. The authors also
thank Irine Peng for comments on an early version of the
manuscript.}

\section{Results from \cite{EFW1} and what remains to be done}
\label{section:suffices:to:prove}

\noindent{\bf Remark:} All terminology in the following theorems
is defined in \cite{EFW1}.  Most of it is recalled in
\S\ref{section:prelim} below. In particular, whenever we wish to
make a statement that refers to either $\Sol$ or $\DL(m,m)$ we
will use the notation $X(m)$ and refer to the space as the {\em
model space}. As in \cite{EFW1}, $\Sol(m)$ denotes $\Sol$ with the
dilated metric
\begin{displaymath}
ds^2 = dz^2 + e^{-2mz} dx^2 + e^{2mz} dy^2.
\end{displaymath}

The main result of this paper is the following.  The analogue of
this theorem for $X(m,n)$ is proved in \cite[Section 5]{EFW1}.

\begin{theorem}
\label{theorem:main:step:two} For every $\delta > 0$, $\kappa > 1$
and $C > 0$ there exists a constant $L_0 > 0$ (depending on
$\delta$, $\kappa$, $C$) such that the following holds: Suppose
$\phi: X(n) \to X(n')$ is a $(\kappa,C)$ quasi-isometry. Then for every
$L > L_0$ and every box $B(L)$, there exists a subset $U \subset
B(L)$ with $|U| \ge (1 - \delta) |B(L)|$ and a height-respecting map
$\hat{\phi}(x,y,z) = (\psi(x,y,z),q(z))$ such that
\begin{itemize}
\item[(i)]
\begin{displaymath}
d(\phi|_U,  \hat{\phi}) = O(\delta L).
\end{displaymath}
\item[(ii)] For $z_1, z_2$ heights of two points in $B(L)$, we
have
\begin{equation}
\label{eq:q:bilishitz} \frac{1}{2\kappa} |z_1 - z_2| - O(\delta L) <
|q(z_1) - q(z_2)| \le 2 \kappa |z_1 - z_2| + O(\delta L).
\end{equation}
\item[(iii)] For all $x \in U$, at least $(1-\delta)$ fraction of
the vertical geodesics passing within $O(1)$ of $x$ are $(\eta,
O(\delta L))$-weakly monotone.
\end{itemize}
\end{theorem}

This theorem, combined with results in \cite[Section 6]{EFW1}
proves that any quasi-isometry $\phi:X(m){\rightarrow}X(m')$ is
within bounded distance of a height respecting quasi-isometry.
This is done in two steps there, the first stated as \cite[Theorem
6.1]{EFW1} roughly shows that $\phi$ respects height difference to
sublinear error. Then in \cite[Section 6.2]{EFW1} we give an
argument that shows this implies $\phi$ is at bounded distance
from height respecting.  The deduction of Theorem
\ref{theorem:sol} from this fact is already given explicitly in
\cite[Section 7]{EFW1}.

The proof of Theorem \ref{theorem:main:step:two} uses the
following consequence of \cite[Theorem 4.1]{EFW1}:

\begin{theorem}
\label{theorem:partIinImage} Suppose $\epsilon, \theta > 0$. Let
$\phi: \Sol \to \Sol$ be a $(\kappa,C)$ quasi-isometry. Then for any
$L'$ sufficiently large  (depending on $\kappa$, $C$, $\theta$),
there exists constants  $R$ and $L$ with $C \ll R \ll L \ll L'$ and
$e^{\epsilon R} \GG L'$ such that for any box $B(L')$ there exist a
collection of disjoint boxes $\{B_i(R)\}_{i \in I}$, a subset $I_g$
of $I$, and for each $i \in I_g$ a subset $U_i \subset B_i(R)$ with
$|U_i| \ge (1-\theta) |B_i(R)|$ such that the following hold:
\begin{itemize}
\item[{\rm (i)}]
\begin{displaymath}
\bigsqcup_{i \in I_g} B_i(R) \subset \phi^{-1}(B(L')) \subset
\bigsqcup_{i \in I} B_i(R),
\end{displaymath}
\item[{\rm (ii)}]
\begin{displaymath}
| \bigsqcup_{i \in I_g} U_i | \ge (1 - \theta) |\phi^{-1}(B(L'))|
\text{
  and } \left| \phi\left(\bigsqcup_{i \in I_g} U_i \right) \right| \ge (1 - \theta) |B(L')|
\end{displaymath}
\item[{\rm (iii)}] For each $i \in I_g$ there exists a product map
$\hat{\phi}_i: B_i(R) \to \Sol$  such that
\begin{displaymath}
d(\phi|_{U_i}, \hat{\phi}_i) = O(\epsilon R).
\end{displaymath}
\end{itemize}
\end{theorem}

\bold{Proof.} Choose $L$ large enough that \cite[Theorem
4.1]{EFW1} holds with the given $\epsilon$ and some $\theta_0 <
\theta$ for any box of size $L$. We cover $\phi{\inv}(B(L'))$ by
boxes of size $L$ in the domain.  Because $\phi$ is a
quasi-isometry, $\phi{\inv}(B(L'))$ is a F\"{o}lner set which
allows us to cover $\phi{\inv}(B(L'))$ by $\cup_{k \in K}B_k(L)$
such that the measure of  $\cup_{k \in K}B_k(L) -
\phi{\inv}(B(L'))$ is small provided $L' \GG L$.  We apply
\cite[Theorem 4.1]{EFW1} to the finite family of boxes $\{B_k(L)|k
\in L\}$ and let $I_g$ be the good boxes which we index without
reference to $k$. By choosing $\theta_0$ small enough and using
the F\"{o}lner condition on $\phi{\inv}(B(L'))$, it is easy to see
that the conclusions of the theorem are satisfied. \qed\medskip

\bold{Recommendations to the reader:}  We strongly recommend that
the reader study \cite{EFW1} before this paper.  In reading this
paper, we recommend that the reader first assume that the map
$\phi$ restricted to each $U_i \subset B_i(R)$ for $i{\in}I_g$ is
within $O(\epsilon R)$ of $b$-standard map, or better yet, the
identity. (Replacing a $b$-standard map with the identity amounts
to composing with a quasi-isometry of controlled constants and so
has no real effect on our arguments.)  This allows the reader to
become familiar with the general outline of our arguments without
becoming too caught up on technical issues.

The reader familiar with \cite{EFW1} can then read
\S\ref{section:prelim} and essentially all of \S
\ref{section:thework}, skipping \S\ref{sec:measure:bilipshitz}
entirely.  In first reading \S\ref{section:prelim}, the reader
might initially read \S \ref{subsection:recollections} through
\S\ref{subsec:trapping} and skip \S\ref{subsection:tangle}. This
last subsection is only required in the case of solvable groups
and then only at the very end of \S\ref{section:families}.  As
remarked there, some of the definitions in
\S\ref{subsection:shadows} may also be omitted on first reading.

\bold{Remarks on the proof:} It is possible to rewrite the
arguments here and first prove that $\phi$ restricted to $U_i
\subset B_i(R)$ for $i{\in}I_g$ is within $O(\epsilon  R)$ of a
$b$-standard map. However, the arguments needed to prove this,
while  not so different in flavor from the arguments in
\S\ref{sec:measure:bilipshitz}, are extremely intricate and
technical.  The proof given here, while slightly more difficult in
some later arguments, is essentially the same proof one would give
after proving that fact.  See \S\ref{section:thework} for more
discussion.

\section{Geometric preliminaries}
\label{section:prelim}

In this section, we describe some key elements of the spaces we
consider.  There is some duplication with \cite{EFW1}, but the
emphasis here is different.

\subsection{Boxes, product maps, almost product maps}
\label{subsection:recollections}

We recall the notion of a box from \cite{EFW1}, first in
$\Sol(m)$. Let
$$ B(L,\vec{0}) =
[-\frac{e^{2mL}}{2},\frac{e^{2nL}}{2}] \times
[-\frac{e^{2mL}}{2},\frac{e^{2nL}}{2}] \times
[-\frac{L}{2},\frac{L}{2}].$$
In our current setting,
$|B(L,\vec{0})| \approx L e^{2L}$ and $Area(\partial B(L,\vec{0}))
\approx e^{2L}$, so $B(L)$ is a F\"olner set.

To define the analogous object in $\DL(m,m)$, we look at the set
of points in $\DL(m,m)$ we fix a basepoint $(\vec{0})$ and a
height function $h$ with $h(\vec{0})=0$. Let $L$ be an even
integer and let $\DL(m,m)_L$ be the
$h{\inv}([-\frac{L+1}{2},\frac{L+1}{2}])$. Then $B(L,\vec{0)})$ is
the connected component of $\vec{0}$ in $\DL(m,m)_L$.  We are
assuming that the top and bottom of the box are midpoints of
edges, to guarantee that they have zero measure.

We call $B(L,\vec{0})$ a box of size $L$
centered at the identity. In $\Sol$, we define the box of size $L$
centered at a point $p$ by $B(L,p)=T_pB(L,\vec{0})$ where $T_p$ is
left translation by $p$. We frequently omit the center of a box in
our notation and write $B(L)$.  For the case of $\DL(m,m)$ it is
easiest to define the box $B(L,p)$ directly. That is let
$\DL(m,m)_{[h(p)-\frac{L+1}{2},h(p)+\frac{L+1}{2}]}=h^{\inv}([h(p)-\frac{L+1}{2},h(p)+\frac{L+1}{2}])$
and let $B(L,p)$ be the connected component of $p$ in
$\DL(m,m)_{[h(p)-{\frac{L+1}{2}},h(p)+\frac{L+1}{2}]}$.  It is
easy to see that isometries of $\DL(m,m)$ carry boxes to boxes.

For $\Sol$, we write $B(R) = S_X \cross S_Y \cross S_Z$. We think
of $S_X$ as a subset of the lower boundary, $S_Y$ as a subset of
the upper boundary, and $S_Z$ as a subset of $\reals$. In the
$DL(n,n)$ case, by $S_X \cross S_Y \cross S_Z$ we mean the set
$\{p \in DL(n,n) \st h(p) \in S_Z \}$ intersected with the union
of all vertical geodesics connecting points of $S_X$ to points of
$S_Y$. We also write $S_Z = [h_{bot},h_{top}]$.  We will use the
notation $\partial^+X$ for the upper boundary and $\partial_{-}X$
for the lower boundary.

\begin{definition}[Product Map, Standard Map]
A map $\hat{\phi}: \Sol \to \pSol$ is called a {\em product map}
if it is of the form
$$(x,y,z) \to (f(x),g(y),q(z)) \qquad\text{ or } \qquad (x,y,z)
\to (g(y),f(x),q(z)),$$
where $f$, $g$ and $q$ are functions from
$\reals \to \reals$. A product map $\hat{\phi}$ is called {\em
$b$-standard} if it is the compostion of an isometry with a map of
the form $(x,y,z) \to (f(x),g(y),z)$, where $f$ and $g$ is
Bilipshitz with the Bilipshitz constant bounded by $b$.
\end{definition}

The discussion of standard and product maps in the setting of
$\DL(m,m)$ is slightly more complicated. We let $\Qa_m$ be the
$m$-adic rationals.  The complement of a point in the boundary at
infinity of  $T_{m+1}$ is easily seen to be $\Qa_m$.  Let $x$ be a
point in $\Qa_m$ viewed as the lower boundary, and $y$ a point in
$\Qa_m$ (viewed as the upper boundary). There is a unique
vertical geodesic in $\DL(m,m)$ connecting $x$ to $y$.  To specify
a point in $\DL(m,m)$ it suffices to specify $x,y$ and a height
$z$.  We will frequently abuse notation by referring to the
$(x,y,z)$ coordinate of a point in $\DL(m,m)$ even though this
representation is highly non-unique.

We need to define product and standard maps as in the case of
solvable groups, but there is an additional difficulty introduced
by the non-uniqueness of our coordinates. This is that maps of the
form $(x,y,z) \to (f(x),g(y),q(z))$, even when one assumes they
are quasi-isometries, are not well-defined, different coordinates
for the same points will give rise to different images. We will
say a quasi-isometry $\psi$ is {\em at bounded distance} from a
map of the form $(x,y,z) \to (f(x),g(y),q(z))$ if $d(\psi(p),
(f(x), g(y), q(z)))$ is uniformly bounded for all points and all
choices $p=(x,y,z)$ of coordinates representing each point.  It is
easy to check that $(x,y,z) \to (f(x),g(y),q(z))$ is defined up to
bounded distance if we assume that the resulting map is a
quasi-isometry.  The bound depends on $\kappa,C,m,n,m'$ and $n'$.

\begin{definition}[Product Map, Standard Map]
A map $\hat{\phi}: \DL(m,m) \to \DL(m',m')$ is called a {\em
product map} if it is within bounded distance of the form $(x,y,z)
\to (f(x),g(y),q(z))$ or or $(x,y,z) \to (g(y),f(x),q(z))$, where
$f: \ratls_m \to \ratls_{m'}$ (or $\ratls_{n'}$), $g: \ratls_n \to
\ratls_{n'}$ (or $\ratls m'$) and $q : \reals \to \reals$.
 A product map $\hat{\phi}$ is called {\em
$b$-standard} if it is the compostion of an isometry with a map
within bounded distance of one of the form $(x,y,z) \to
(f(x),g(y),z)$, where $f$ and $g$ are Bilipshitz with the
Bilipshitz constant bounded by $b$.
\end{definition}

\begin{definition}
\label{definition:almost:product:map} Given a quasi-isometric
embedding $\phi:B(R){\rightarrow}X(n')$, we say $\phi$ is an
$(\alpha,\theta)$ almost a product map if  there exist subsets $U
\subset B(R), E_1 \subset S_X$ and $E_2 \subset S_Y$ of relative
measure $1-\theta$ such that $U = \{ (x,y,z) \st x \in E_1, y \in
E_2, z \in S_Z \}$ and all geodesics connecting points in $E_1$ to
points in $E_2$ have $\epsilon$ monotone images under $\phi$.
\end{definition}

\bold{Remark.} We think of $f$ and $g$ as defined only on $E_i$.
So by $f(I)$ we mean $f(I \cap E_1)$.

\begin{lemma}
\label{lemma:almost:product:map} Given a $(\alpha,\theta)$-almost
product map $\phi$ there exists a subset $U^* \subset U$ with
relative measure $1-128\theta^{\frac{1}{2}}$ and a (partially
defined) product map $\hat \phi:U^* {\rightarrow} X(m')$ such that

$$d(\phi|_U(p),\hat \phi(p)) \leq \alpha R$$

\noindent for all $p{\in}U$.
\end{lemma}

\bold{Proof.} This is the content of \cite[Lemma 4.11 and
Proposition 4.12]{EFW1} \qed

\bold{Remark.} With an appropriate choice of constants, the
converse of Lemma \ref{lemma:almost:product:map} is also true.

\bold{Notation.} Using Lemma \ref{lemma:almost:product:map}, we
write an (almost) product map $\hat{\phi}: B(R) \subset X(n) \to
X(n')$ as $(x,y,z) \to (f(x),g(y),q(z))$, so the domain of $f$ is
$S_X$ etc. We will always work with (almost) product maps of this
kind, the arguments for those of the form $(x,y,z) \to (f(y),
g(x), q(z))$ are almost identical.  One can also formally deduce
any result we need about almost product maps of the form $(x,y,z)
\to (f(y), g(x), q(z))$ from the analogous fact about those of the
form $(x,y,z) \to (f(x),g(y),q(z))$ by noting that these two forms
of almost product map differ by either pre- or post-composition
with an isometry.

\subsection{Discretizing $\Sol$}
\label{subsection:discretization}

In this subsection, we introduce a discrete model for $\Sol(n)$
which has some technical advantages at some points in the
argument.  We will often make arguments for the discrete model
instead of for $\Sol(n)$ itself. The discrete model is
quasi-isometric to $\Sol(n)$ and in fact $(1,\rho_1)$
quasi-isometric for a parameter $\rho_1$ which will we choose so
that $C \ll \rho_1 \ll \epsilon R$.

The basic idea is to take a $\rho_1$ net in $\Sol(n)$ and replace
$\Sol(n)$ by a graph on this net.  It is possible to do this by
taking a arithmetic lattice in $\Sol$, taking a deep enough
congruence subgroup, and taking the Cayley graph.  More
concretely, we write $\Sol(n)$ as $\Ra{\ltimes}\Ra^2$, and
consider $\rho_1 \Za \subset \Ra$ and $\rho_1 \Za^2 \subset
\Ra^2$.  Here we view $\Ra^2 \subset \Sol(n)$ as the plane at
height zero.  We then form a $\rho_1$ net in $\Sol(n)$ by taking
the union
$$\mathcal{G}=\bigcup_{a \in \rho_1 \Za} a{\cdot}{\rho_1 \Za^2}.$$

To make this a graph, we connect by an edge any pair of points in
$\mathcal{G}$ whose heights differ by $\rho_1$ and which are
within $10\rho_1$ of one another.  We metrize this graph by
letting lengths of edges be the distance between the corresponding
points in $\Sol(n)$, so all edges have length $O(\rho_1)$.

We can also replace $DL(m,m)$ with a graph whose edges have length
$O(\rho_1)$. For this we assume $\rho_1{\in}\Na$.  Consider only
vertices in $DL(m,m)$ in $h{\inv}(\rho_1\Za)$.  Join two vertices
by an edge of length $\rho_1$ if there is a monotone vertical path
between them. The resulting graph is clearly quasi-isometric to
$DL(m,m)$ and is in fact $DL(m^{\rho_1},m^{\rho_1})$ but with the
edge length fixed as $\rho_1$ instead of $1$.

We remark here that constants that are said to depend only on $K,C$
and the model geometries often also depend on the discretization
scale. This is because the discretization process effectively
replaces the model space with a graph.

\subsection{Shadows, slabs and coarsenings}
\label{subsection:shadows}

\bold{Shadows and projections:} Let $H$ be a subset of an
$y$-horocycle, and suppose $\rho > 1$. By the $\rho$-shadow of
$H$, denoted $\Shadow(H,\rho)$, we mean the union of the vertical
geodesic rays which start within distance $\rho$ of $H$ and go
down. If $H$ is a $x$-horocycle, then the we use the same
definition except that the geodesic rays are going up.

Given a subset of a $y$-horocycle $H$, we let
$\pi_-(H)=\partial_-X{\cap}\Sh(H,\rho_1)$.  We define $\pi^+(H)$
for a subset of an $x$-horocycle similarly. Note that we are
suppressing $\rho_1$ in the notation.  In any context where
$\pi^+$ or $\pi_-$ are used, $\rho_1$ will be fixed in advance.

\bold{The number $\Delta(H)$.}
For a horocycle $H$ in a box $B(R)$, let $\Delta(H) = \min(h_{top}
- h(H), h(H) - h_{bot})$. Thus, $\Delta(H)$ measures how far is $H
\cap B(R)$ is from the  top and bottom  of $B(R)$.

\bold{The branching numbers $B_X$ and $B_X'$.}
We define $B_X$ to be the branching
constant of $X$. For solvable Lie groups $B_{X(n)}={n}$, for
Diestal-Leader graphs, $B_{X(n)}={\log(n)}$. We use the shorthand
$B_X'$ for $B_{X(n')}$.

\bold{Measures on the boundary at infinity.} Note that the boundaries
$\partial_-(X)$ and $\partial^+(X)$ are homogeneous spaces, and thus
have a natural Haar measure. (This measure is Lebesque measure on
$\reals$ if $X = \Sol$ and the natural measure on the Cantor set if $X
= DL(n,n)$). We normalize the measures
by requiring that the shadow of a point at height $0$ has measure $1$.
These measures are all denoted by the symbol $| \cdot |$.
Note that for any point $p \in X$,
\begin{equation}
\label{eq:measure:shadow:point} |\pi_-(\{p\})|e^{-B_X h(p)}=1
\end{equation}

\bold{The parameter $\beta''$.} We choose an arbitrary $\beta''$
with $\beta'' \ll
1$, with the understanding that $\epsilon$ and $\theta$ will be chosen
so that $\epsilon \ll \beta''$ and $\theta \ll
\beta''$. The parameter $\beta''$ will be fixed until
\S\ref{sec:subsection:end:proof:theorem2.1}.

\bold{Slabs:} The objects we refer to as {\em slabs} will always
be subsets of the part of the  box $B(R)$ which is at least
$4\kappa^2\beta''R$ from the boundary of $B(R)$, will always be
defined in reference to a horocycle $H$ in $B(R)$, and are always
contained in $Sh(H,\rho)$.  We give definitions only for slabs in
shadows of $y$ horocycles, those for $x$ horocycles are analogous
and can be obtained by applying an appropriate flip.  If we choose
$h_2 < h_1 < h(H)$, a {\em slab }in $B(R)$ below $H$ is the subset
$Sl_{2}^{1}(H)$ which is defined to be the subset of the shadow of
$H$ which is between heights $h_2$ and $h_1$.

\bold{Recommendation to the reader:} The remainder of this
subsection might be omitted on first reading.

\bold{Coarsening:} In order to work with more regular sets, we
define an operation to {\em coarsen} subsets of either boundary.

Let $a_1,a_2$ be two points in a (log model) hyperbolic plane
(which we think of as the $xz$ plane in $\Sol$). Let
$h^+(a_1,a_2)$ be the  height at which vertical geodesics leaving
$a_1$ and $a_2$ are one unit apart.  This function clearly extends
to the lower boundary of the hyperbolic plane.  We further extend
the function to $\Sol$ by letting
$h^+(p_1,p_2)=h^+(\pi_{xz}(p_1),\pi_{xz}(p_2))$ where
$\pi_{xz}(x,y,z)=(x,z)$. If $I = [a,b]$ is an interval, we write
$h^+(I)$ for $h^+(a,b)$. Note that we can define $h_{-}$ similarly
in a $yz$ plane. For $\DL(n,n)$ we define $h^+(a_1,a_2)$ as the
height in $T_n$ at which vertical geodesics leaving $a_1$ and
$a_2$ meet.  Again $h_-$ is defined similarly.

The operation of coarsening replaces any set $F$ by a set
$\cC_z(F)$ which is a union of open intervals  of a certain size
depending on $z$. For $F \subset
\partial_-X$ and $z \in \reals$, let $\cC_z(F)$ denote the set
of $x \in \partial_-X$ such that there exists $x' \in F$ with
$h^+(x,x') < z$. Similarly, for $F \subset \partial^+X$ and $z \in
\reals$, let $\cC_z(F)$ denote the set of $y \in
\partial^+(X)$ such that there exists $y' \in F$ with $h^-(y,y') >
z$.

\bold{Generalized Slabs:} Given two sets $E^+{\subset}\partial^+X$
and $E_-{\subset}\partial_-X$, and two heights $h_2<h_1$, we
define a set
$$S(E_{-}, E^+, h_2, h_1) = \{ (x,y,z)  \text{ such that }
h_2<z<h_1  \text{ and } x {\in} E_{-}, y {\in} E_{+}\}.$$

In words $S(E_-, E^+, h_2, h_1)$ is the set of points on geodesics
joining $E^+$ to $E_{-}$ with height between $h_1$ and $h_2$.  We
refer to these sets as {\em generalized slabs}, though in general
there geometry can be very bad, depending on the geometry of
$E^{+}$ and $E_{-}$.  Generalized slabs will always be subsets of
the part of the  box $B(R)$ which is at least $4\kappa^2\beta''R$
from the boundary of $B(R)$, even if this is not explicit in our
specification of $E^+$ and $E^-$. In particular, slabs as defined
above are special cases of generalized slabs, with $Sl_2^1(H) =
S(\pi_-(H), S_Y , h_2, h_1)$ where $h_2<h_1<h(H)$.

Clearly boxes are very special generalized slabs, and
we prefer to work in general with generalized slabs that are unions
of boxes. One can obtain a generalized slab that is a union of boxes
by coarsening $E^{+}$ and $E_{-}$.
Let $h_3$ and $h_4$ be two additional heights, and consider $S(
\cC_{h_3}(E_{-}), \cC_{h_4}(E^+), h_2, h_1)$.  Observe that as
long as $h_3 \leq h_2$  and $h_4 \geq h_1$ we have

$$S(E_{-}, E^+, h_2, h_1) = S( \cC_{h_3}(E_{-}), \cC_{h_4}(E^+),
h_2, h_1).$$

We will need some information concerning the geometry of coarse
enough generalized slabs.
\begin{lemma}
\label{lemma:generalized:slabs} Choose $h_3 \geq h_1$ and $h_4
\leq h_2$. Then any generalized slab of the form
$S( \cC_{h_3}(E_{-}),
\cC_{h_4}(E^+), h_2, h_1)$ is a union of boxes of size $h_1-h_2$.
In the $\DL(m,m)$ case, $S( \cC_{h_3}(E_{-}), \cC_{h_4}(E^+), h_2,
h_1)$ is a disjoint union of boxes of size $h_1-h_2$.  In the
$\Sol$ case, $S( \cC_{h_3}(E_{-}), \cC_{h_4}(E^+), h_2, h_1)$ is
not a disjoint union of boxes, but any such set contains a
disjoint union of boxes of height $h_1-h_2$ that contain
$\frac{1}{25}$ of the volume of $S( \cC_{h_3}(E_{-}),
\cC_{h_4}(E^+), h_2, h_1)$.  Furthermore, the number of vertical
geodesics in $S( \cC_{h_3}(E_{-}), \cC_{h_4}(E^+), h_2, h_1)$ is
comparable to:

$$\frac{\Vol(( \cC_{h_3}(E_{-}),\cC_{h_4}(E^+), h_2, h_1))}{h_1-h_2}
e^{B_X(h_1-h_2)}$$

\noindent i.e. it is comparable to the area of the cross-section
times $e^{B_X(h_1-h_2)}$
\end{lemma}

\bold{Proof.} That $S( \cC_{h_3}(E_{-}), \cC_{h_4}(E^+), h_2,
h_1)$ is a union of boxes is clear from the definition of
coarsening. In the $\DL(m,m)$ case, the set between $h_1$ and
$h_2$ is a disjoint union of boxes of size $h_1-h_2$, so the
result follows.  For $\Sol$ one proves the result by considering
the set $W=S( \cC_{h_3}(E_{-}), \cC_{h_4}(E^+), h_2, h_1) \cap
h{\inv}(z)$ for any $z{\in}(h_2,h_1)$. It is clear that $W$ is
covered by it's intersection with boxes of size $h_1-h_2$, all of
which are rectangles of the same size and shape.  Using the Vitali
covering lemma, one finds a subset of the boxes which cover a
fixed fraction of the measure of $W$.  Since the volume of $S(
\cC_{h_3}(E_{-}), \cC_{h_4}(E^+), h_2, h_1)$ is the area of the
cross section times $h_1-h_2$, we are done.

The claim concerning numbers of vertical geodesics is obvious for
a box.  The proof in general can be reduced to that case using the
earlier parts of this lemma.\qed

\subsection{The trapping lemma}
\label{subsec:trapping}

In this subsection we state some results relating to areas,
lengths and shadows.  These are used in the proof of Theorem
\ref{theorem:fullbox:align}.  Some similar statements are
contained in \cite[\S5.2]{EFW1}.

\begin{lemma}
\label{lemma:shadows:coarse:length} Let $Q$ be a subset of an
x-horocycle $H$. Then $\pi_-(Q) = \pi_-(H)$ and
\begin{displaymath}
\ell(Q) \approx |\pi_-(H)| |\pi_+(Q)|,
\end{displaymath}
where by $\ell(Q)$ we mean the length of the intersection of the
$3 \rho$ neighborhood of $Q$ with $H$, and the implied constants
depend on $\rho$.
\end{lemma}

\bold{Proof.} This follows from (\ref{eq:measure:shadow:point}).
\qed\medskip

\begin{lemma}
\label{lemma:length:area} Suppose $\gamma \subset B(R)$ is a path.
Let $L$ be a euclidean plane intersecting $B(R)$, and suppose $U
\subset L \cap B(R)$. Suppose also that any vertical geodesic
segment from the bottom of $B(R)$ to the top of $B(R)$ which
intersects $U$ also intersects the $\rho$-neighborhood of
$\gamma$. Then,
\begin{displaymath}
\ell(\gamma) \ge c(\rho) \Area(U)
\end{displaymath}
(in the above, $c(\rho)$ is a constant, and both the length and
the area are measured using the $X(n)$ metric).
\end{lemma}

\bold{Proof.} First note that if $L'$ is another Euclidean plane,
and $U'$ is the (vertical) projection of $U$ on $L'$, then
$\Area(U') = \Area(U)$.

Now subdivide $\gamma$ into $k$ segments of length $\rho$. Let
$x_i$ be the midpoints of such a segment. Let $R_i$ be a rectangle
at the same height as $i$, such that $x_i$ is in the center of
$R_i$, and the sides of the rectangle have length $2\rho$. Then
the $X(n)$-area of $R_i$ is independent of $i$, and the projection
of the union of the $R_i$ to $L$ must cover $U$. Therefore $k >
c_2(\rho) \Area(U)$, and hence $\ell(\gamma) > c_1(\rho)
\Area(U)$. \qed\medskip

\begin{lemma}[Trapping Lemma]
\label{lemma:trapping} Suppose $\rho_1 \GG C$, $\phi: X(n) \to
X(n)$ is a $(\kappa,C)$ quasi-isometry, and and $H$ is a subset
(not necessarily connected) of an $x$-horocycle in $X(n)$.

Suppose $Q$ is a subset of a finite union of horocycles in $X(n)$,
such that such that the $\kappa\rho_1$-neighborhood of $\phi(Q)$
intersects every vertical geodesic starting from the
$\rho_1$-neighborhood of $H$ and going down. Then,
\begin{displaymath}
\ell(Q) \ge c_1 \ell(H)
\end{displaymath}
where $c_1 = c_1(\rho_1)$.
\end{lemma}

\bold{Proof.} Discretize $H$ on the scale $\rho_1$, and apply
Lemma~\ref{lemma:length:area}. \qed\medskip

Lemma~\ref{lemma:trapping} is sufficient for applications to
$DL(n,n)$. For applications to $\Sol$, we will need a
generalization that is stated in the next subsection.

\subsection{Tangling and generalized trapping}
\label{subsection:tangle}

% \begin{lemma}
% \label{lemma:shadows:distance} Suppose $p \in \Sh(H,\rho)$, and $q$
% is connected to $p$ by a path $\gamma$ of length at most $????
% d(p,q)$ which does not intersect the $2\rho$-neighborhood of $H$.
% Then $q \in \Sh(H,2\rho)$. In particular, if $d(q,p) < d(p,H) - 2
% \rho$, $q \in \Sh(H,2\rho)$.
% \end{lemma}

The following (obvious) result about $DL$ graphs is used
implicitly in the proof of Theorem \ref{theorem:fullbox:align}.

\begin{lemma}
\label{lemma:shadows:disconnect} Suppose $\rho > 1$ and
$p$ and $q$ are two points
in $DL(n,n)$. Suppose also $p \in \Sh(H,\rho)$, $q \in
\Sh(H,\rho)^c$. Then any path connecting $p$ to $q$ passes within
$\rho$ of $H$.
\end{lemma}

{\bf Proof:} The point is simply that if $\pi_T$ is the projection
to the tree $T_{n+1}$ transverse to $H$, then $\pi_T(\Sh(H,\rho))$
is exactly the set directly below the unique point $x$ which is
$\rho$ units above the projection of $\pi_T(H)$.  And removing $x$
disconnects $T_{n+1}$.
\qed\medskip

The lemma above is false for the case of $\Sol$. We will need the
following variant: Fix an integer $\rho >  100$ for the remainder
of this section.

\begin{definition}[Tangle]
\label{def:tangle} Let $H$ be a horocycles. We say that a path
$\bar{\gamma}$  {\em tangles with $H$ within distance $D$} if
either $\bar{\gamma}$ intersects the $\rho$ neighborhood of $H$ or
\begin{displaymath}
\tau(\bar \gamma, H)=\sum_{j=1}^{\frac{D}{\rho}}
\frac{\ell(\bar{\gamma} \cap \{ p \st ja \le d(p,H) \le
  (j+1)3 \rho \})}{\nu(ja)} > 100.
\end{displaymath}
Here $\nu(r)$ is the volume of the ball of radius $r$ in the
hyperbolic plane.  Informally, in order for the path $\bar{\gamma}$
to tangle with the horocycle $H$,
$\bar{\gamma}$ has to spend a lot of time near $H$, with
the closest approaches to $H$ carrying more weight. 

We say $\bar{\gamma}$ {\em tangles with a finite union  of
horocycles $\mathcal H$ within distance $D$} if
\begin{displaymath} \sum_{H{\in}{\mathcal
H}}\tau(\bar{\gamma},H)>100
\end{displaymath}
\noindent where $D$ is implicit in our definition of $\tau$.
\end{definition}

We first state an easy lemma to illustrate situations in which
paths can be forced to tangle with a horocycle.

\begin{lemma}
\label{lemma:shadows:tangle:example} Let $\rho$ be as above and
let $H$ be horocycle in $\Sol$. Suppose $p$ and $q$ are two points
in $\Sol$ such that $p \in \Sh(H,\rho/3)$ and $q \in
\Sh(H,\rho)^c$. Then any path from $p$ to $q$ of length less than
$L$ tangles with $H$ at distance $\log(L)$.
\end{lemma}

\bold{Proof:} This is an easy hyperbolic geometry argument applied
to the projection of the path a hyperbolic plane transverse to
$H$. \qed\medskip

For our applications, we require a more technical variant of
Lemma~\ref{lemma:shadows:tangle:example}.  In our arguments, we
deal with $\Sh(H,\rho_1)$ where $\rho_1$ is the discretization
scale.  For this reason, $\Sh(H,\frac{\rho_1}{3})$ is not a good
notion and we need to specify the set we consider differently.
Given an horocycle $H$ and constant $D'$, we say a point $p$ is
{\em $D'$-deep in $\Sh(H,\rho)$} if $p$ is more than $D'$ below
$H$ and more than $\frac{D}{9}$ from the edges of the shadow.

\begin{lemma}
\label{lemma:shadows:tangle} Let $\rho$ be as above, and choose
constants $\rho \ll D_1 \ll D_2$. Let $H$ be a horocycle and
suppose $p$ and $q$ are two points in $\Sol$ such that $p$ is
$D_2$-deep in $\Sh(H,\rho)$ and $q \in \Sh(H,\rho)^c$. Then any
path from $p$ to $q$ of length less than $e^{D_1}$ tangles with
$\mathcal H$ within distance $D_2$.
\end{lemma}

\bold{Proof:} This is an easy hyperbolic geometry argument applied
to the projection of the path a hyperbolic plane transverse to
$H$. \qed\medskip

For a family $\cF$ of vertical geodesic segments, we let $\|\cF\|$
denote the area of $\cF \cap P$, where $P$ is a Euclidean plane
intersecting all the segments in $\cF$. (If there is no such plane
we break up $\cF$ into disjoint subsets $\cF_i$ for which such
planes exist, and define $\|\cF\| = \sum_i \|\cF_i\|$).

\begin{lemma}[Generalized Trapping Lemma]
\label{lemma:generalized:trapping} Let $\rho \ll D_2$ be constants
as above.  Suppose $\cF$ is a family of vertical geodesic
segments,  suppose $Q$ is a subset of a finite union $\mathcal H$
of horocycles. Suppose also that for each $\gamma \in \cF$,
$\gamma$ tangles with $\mathcal H$ within distance $D_2$ and that
$\gamma$ is contained in $N(Q^c, D_2)$. Then, $\ell(Q) \ge \omega
\|\cF\|$, where $\omega$ depends only on $\kappa,C,n$ and the
constants in the definition of tangle.
\end{lemma}

\bold{Proof.} We assume $P$ is a Euclidean plane intersecting all
the geodesics in $\cF$, the general case is not much harder. Let
$S(r) = \{ p \st r \le d(p,Q) \le
  r+a \}$. Then, $|S(r)| = c \, \ell(Q) \nu(r)$, where $c$ depends only
  on $a$. Then, we have by \cite[Proposition 5.4]{EFW1},
\begin{displaymath}
\ell(Q) = \frac{|S(r)|}{c \nu(r)} \ge \omega_1
\frac{|(S(r))|}{\nu(r)} \ge \omega_1 \int_{\cF\cap P}
\frac{\ell(\gamma \cap (S(r)))}{\nu(r)} \, d\gamma \ge \omega_2
\int_{\cF \cap P}\frac{\ell(\gamma) \cap
  S(r))}{\nu(r)} \, d\gamma
\end{displaymath}
where we have identified the space of vertical geodesics with $P$
and $\omega_1$ and $\omega_2$ depend only on $(\kappa, C, a)$.
After writing $r = ja$, summing the above equation over $j$ and
using the assumption that $\gamma$ tangles with $\mathcal H$
within distance $D_2$ and is contained in $N(Q^c,D_2)$ for all
$\gamma \in \cF$, we obtain that $\ell(Q) \ge \omega |\cF \cap P|$
as required. \qed\medskip

\section{Improving almost product maps}
\label{sec:measure:bilipshitz}

In this section, we make some arguments that improve the
information available concerning $\phi|_{B_i(R)}$ where
$i{\in}I_g$.  More or less, by throwing away another set of small
measure, we show that $\phi$ maps many slabs to particular nice
generalized slabs.  We also show that the map $q$ can be taken to
be a linear map.

\bold{Recommendation to the reader:} The reader may wish to skip
this section on first reading and continue reading assuming that
$\phi|_{U_i}$ is $b$-standard or within $O(\epsilon R)$ of  a
$b$-standard map. All the results in this section are somewhat
technical in nature.

\subsection{Bilipschitz in measure bounds}
\label{subsection:measure:bilipschitz}

It is clear that the image of a slab under a product map is a
generalized slab and that the image of a slab under a $b$-standard
map is a slab.  We need to work instead with images of slabs under
almost product maps.  Given an almost product map
$\phi:B(R){\rightarrow}\Sol$ one wants to understand the image of
$Sl_2^1(H)$. In general, there is not an obvious relation between
$\phi({Sl}_{2}^{1}(H))$ and $S(f(\pi_-(H)), g(S_Y) , q(h_2),
q(h_1))$. We will show that this is true, at least up to sets of
small measure, for appropriately chosen slabs, once we coarsen the
image of the slab. To this end we let $h=h(H)$ and fix a height
$h_1<h$ and define:

\begin{equation}
\label{equation:slab:containments} \hat Sl_{2}^{1}(H)=S(
\cC_{q(h_1)}(f(\pi_-(H)), \cC_{g(h_2)}(g(S_Y)), q(h_2), q(h_1)).
\end{equation}

\noindent Note that $\hat Sl_{2}^{1}(H)$ is a union of boxes of
size $q(h_1)-q(h_2)$.   When the choice of $H$ is clear, we
suppress reference to $H$ and consider $h_1,h_2$ and also write
$h$ for $h(H)$.

In this section, we prove two lemmas which show that we can
restrict attention to $Sl_2^1(H)$ which are almost entirely in
$U^*$ and whose (coarsened) image is mostly a collection of boxes
contained in (a small neighborhood) of the image of $U^*$.

\bold{Terminology:} In order to discuss properties of
$\hat{Sl}_2^1(H)$ without fixing either the orientation of $H$ or
the almost product map on $B_i(R)$, we introduce some terminology.
This terminology is justified by comparison with the case where
$\hat{Sl}_2^1(H)$ is a slab.  We refer to the direction in $z$
that goes from $q(h_2(H))$ to $q(h(H))$ as {\em towards the
horocycle} and the opposite direction as {\em away from the
horocycle}. Similarly, there is a direction, either $x$ or $y$
that one can think of as being {\em along the horocycle} where the
other direction is {\em transverse to the horocycle}. If $H$ is an
$x$ horocycle and our product map is of the form
$(x,y,z){\rightarrow}(f(x),g(y),q(z))$ then $x$ is along the
horocycle and $y$ is transverse to the horocycle.

Let $\phi$ be an $(\epsilon, R)$ almost product map and $\hat
\phi$ the corresponding (partially defined) product map. The
following equation follows from the definitions. It says that the
image of the intersection of certain slabs with the good set is
essentially contained in a corresponding slab.
\begin{equation}
\label{eq:phi:tilde:B:subset} \phi(Sl_{2}^{1}(H) \cap U_*) \subset
N_{O(\epsilon R)} \hat{\phi}(Sl_2^1(H) \cap U_*)  \subset
N_{O(\epsilon R)}(\hat Sl_2^1(H)).
\end{equation}

\noindent The following two lemmas yield a strengthening of the
equation above.  The first lemma provides a lower bound on the
measure of $Sl_{2}^{1}(H) \cap U_*$ and so on the measure of
$N_{O(\epsilon R)}(\hat Sl_2^1(H))$ for most choices of $H$. The
second lemma provides an upper bound on $N_{O(\epsilon R)}(\hat
Sl_2^1(H))$ and even $N_{O(\rho_1)}(\hat Sl_2^1(H))$ for a more
restricted set of choices of $H$. To do this, we  actually need to
modify $\hat Sl_2^1(H)$ in a way that we describe in Lemma
\ref{lemma:measure:upper:bilip}.

Given any subset $A \subset B(R)$ and any constant $d<1$, we
denote by $A^d$ the intersection of $A$ with the points in $B(R)$
more than $dR$ of the $\partial B(R)$.

\medskip

\begin{lemma}
\label{lemma:measure:lower:bilip}  Given $\beta' \GG \beta \GG
\alpha \GG 1$, there exist constants $c_1,c_2$ depending on
$\epsilon,\theta$ and $\beta'$ and a subset a subset $E_{**}$ of
$S_X$ with $|S_X\setminus E_{**}| \le c_1(\theta,
\epsilon,\beta')|S_X|$ with the following properties. Given a
$y$-horocycle $H$ intersecting $B(R)$ more than $2 \kappa \beta'
R$ away from $\partial B(R)$ and with $\pi_-(H)$ containing a
point of $E_{**}$ and any slab $Sl_2^1(H)$ such that $\beta' R >
|h_1(H) - h_2(H)| > \beta R$, $4 \beta R > |h(H) - h_1(H)|
>  2\alpha R$, we have
\begin{equation}
\label{eq:tilde:B:H:cap:U:star} |Sl_2^1(H) \cap U_*| \ge (1 -
c_2(\theta,\epsilon)) |Sl_2^1(H)|.
\end{equation}
\end{lemma}

Our current notion of $\hat{Sl}_2^1(H)$ is a bit too coarse. In
particular, there can be points in this set that are $O(R)$ away
from $\phi(SL_2^1(H))$.  We introduce some notations needed to
describe a subset of $\hat{Sl}_2^1(H)$ which can be controlled
more easily. Given a set $D\subset B(R)$, we denote by
$S_Y{\cap}D$ the set in $S_Y$ consisting of $y$ coordinates of
points in $D$. We then define
$$\tilde Sl_{2}^{1}(H,D)=S( \cC_{q(h_c)}(f(\pi_-(H)), \cC_{g(h_2)}(g(S_Y{\cap} D)),
q(h_2), q(h_1)).$$  The fact that we only intersect the $y$
coordinate with $D$ is not an accident, it is due to the fact that
we consider sets which are ``large" in the $y$ direction and
``small" in the $x$ direction.

\begin{lemma}
\label{lemma:measure:upper:bilip} Given $\beta'' \GG \beta' \GG
\beta \GG \alpha \GG 1$, there exist constants $c_3,c_4$ depending
on $\epsilon,\theta$ and $\beta'$ and a subset a subset $E_*$ of
$S_X$ with $|S_X\setminus E_*| \le c_3(\theta, \epsilon)|S_X|$
with the following properties. For any $y$-horocycle $H_0$
intersecting $B(R)$ more than $ 4 \kappa^2 \beta'' R$ away from
$\partial B(R)$ with $\pi_-(H_0)$ containing a point of $E_*$,
consider all horocycles $H$ in $S=S(\pi_-(H_0), S_Y, h(H_0),
h(H_0) + \beta''R){\cap}U_*$ with $\pi_-(H)$ containing
a point of $E_*$ and any constants $\beta' R
> |h_1(H) - h_2(H)| > \beta R$, $4 \beta R > |h(H) - h_1(H)|
>  \alpha R$ such that the slab $Sl_2^1(H)$ is also contained in $S$,
letting $\tilde {Sl}_2^1(H)=\tilde{Sl}_2^1(H, \phi(S))$, we have:
\begin{equation}
\label{eq:hat:B:H:cap:N} |\tilde{Sl}_2^1(H) \cap N_{\rho_1}(\phi(U_*
\cap Sl_2^1(H)))| \ge (1 - c_4(\theta,\epsilon))
|\tilde{Sl}_2^1(H)|.
\end{equation}
For $i=3,4$, we have $c_i(\theta,\epsilon) =c_i(\theta, \epsilon)
\to 0$ as $\epsilon \to 0$ and $\theta\to 0$.
\end{lemma}

Saying $H_0$ intersects $B(R)$ more than $4\beta''^2 R$ from
$\partial B(R)$ is the same as saying $H_0$ intersects the box
$B((1-2\beta'')R)$ with the same center as $B(R)$. The point is to
stay away from the edge of the box. See the remarks in the
definition of slabs and generalized slabs.

While the proof of Lemma \ref{lemma:measure:lower:bilip} is essentially an application of the Vitali covering lemma, the proof of Lemma \ref{lemma:measure:upper:bilip} depends on the fact that quasi-isometries roughly preserve volume.  We will also need this fact to deduce some corollaries from Lemma \ref{lemma:measure:lower:bilip} and \ref{lemma:measure:upper:bilip}.  We recall a precise
statement from \cite{EFW1}:

\begin{proposition}
\label{prop:qivolume} Let $\phi:X{\rightarrow}X'$ be a
continuous $(\kappa,C)$ quasi-isometry.  Then for any $a \GG C$
there exists $\omega_1 > 1$ with $\log \omega_1 = O(a)$ such that
for any $U \subset X$,
\begin{displaymath}
\omega_1^{-1} | \phi(N_a(U)) | \le | N_a(U)| \le \omega_1
|N_a(\phi(U))|
\end{displaymath}
where $N_a(U) = \{ x \in X \st d(x,U) < a \}$.
\end{proposition}

As explained  in \cite{EFW1}, this fact holds much more generally for metric measure
spaces which satisfy relatively mild conditions on the growth of balls.

Before proving the lemma, we state and prove a corollary
concerning measures of cross sections. We note that by the
definitions of the measures on the boundary, for a generalized
slab $S(E_-, E^+, h_2,h_1)$, and $h_1 < z < h_2$, the area (or
equivalently volume) of the $O(1)$ neighborhood of the cross
section at height $z$ (i.e. of $S(E_-, E^+, h_2,h_1) \cap
h^{-1}(z)$) is $|\cC_z(E_+)| |\cC_z(E_-)|$.

%In $\Sol$, we let
%$S_Y'' =[-e^{2(1-\beta'')R},e^{2(1-\beta'')R}]$, and in $\DL(n,n)$
%we let $S_Y'' = S_Y$.   This difference reflects
%the fact that horocycles in $\Sol$ have ends near the edges of
%boxes while horocycles in $\DL(n,n)$ do not.

\begin{corollary}
\label{cor:measure:lower:bilp} Assume $H$ satisfies the hypotheses
of Lemmas \ref{lemma:measure:lower:bilip} and
\ref{lemma:measure:upper:bilip}.  Let $w_1,w_2$ be such that
$2\kappa\alpha R < |h(H)-w_1| < \frac{2}{\kappa}\beta R$ and
$2\kappa\beta R < |w_2-w_1| < \frac{1}{2 \kappa}\beta' R$.  Then:
\begin{equation}
\label{eq:measure:lower:bilip} |\cC_{w_1}(f(\pi_-(H)))|
|\cC_{w_2)}(g(S_Y\cap S)| \ge \omega |\pi_-(H)| |S_Y|,
\end{equation}

and

\begin{equation}
\label{eq:prelim:upper:measure:lip} |\cC_{w_1}(f(\pi_-(H)))|
|\cC_{w_2}(g(S_Y\cap {S}))| \le b |\pi_-(H)||S_Y|,
\end{equation}
where $\omega$ and $b$ depend only on $\kappa$ and $C$.
\end{corollary}

\bold{Proof of Corollary.}
Note that from the structure of $U$ and the fact that $\phi$ is a
quasi-isometry, it follows that for $z_1, z_2 \in S_Z$, we have
\begin{equation}
\label{eq:q:part15:bilishitz} \frac{1}{2\kappa} |z_1 - z_2| -
O(\epsilon R) < |q(z_1) - q(z_2)| \le 2 \kappa |z_1 - z_2| +
O(\epsilon R).
\end{equation}
In particular, $q$ is essentially monotone (up to $O(\epsilon R)$
error).

Given $w_1,w_2$ as in the Corollary, there exist heights $h_1(H)$
and $h_2(H)$ satisfying the hypotheses of Lemmas
\ref{lemma:measure:lower:bilip} and
\ref{lemma:measure:upper:bilip} such that
$q(h_1(H))=w_1,q(h_2(H))=w_2$. We apply those lemmas to the
resulting $Sl_2^1(H)$ and $\tilde{Sl}_2^1(H)$.

Let $\Vol'(X)=|N_{\rho_1}(X)|$ with $\epsilon R \GG \rho_1 \GG C$.
Recall that $|h_1-h_2|>\beta R$ for some $\beta \GG \epsilon$.  By
Lemma~\ref{lemma:generalized:slabs} and the fact that the measure of
the $O(\epsilon R)$ neighborhood of a box of size $\beta R$ is
comparable to the measure of a box of size $\beta R$, we have
\begin{equation}
\label{eq:vol:gen:slab:epsilon:nbhd} (1-c) \Vol'(N_{O(\epsilon
R)}(\tilde{Sl}^1_2(H))) \le \Vol'(\tilde{Sl}^1_2(H)),
\end{equation}
where $c$ is
a constant that depends only on $\frac{\epsilon}{\beta}$ and which
goes to $0$ as $\epsilon$ goes to zero.

Note that  (\ref{eq:phi:tilde:B:subset}) continues to hold when we
replace $\hat{Sl}_2^1(H)$ by $\tilde{Sl}_2^1(H)$. Therefore, by
(\ref{eq:phi:tilde:B:subset}) and (\ref{eq:hat:B:H:cap:N}) we have
\begin{equation}
\label{eq:slab:volume:comparison} (1-c_3) (1-c) \Vol'(\tilde
Sl_2^1(H)) \le (1-c) \Vol'( \phi(Sl_{2}^{1}(H) \cap U_*)) \leq
\Vol'(\tilde Sl_2^1(H)),
\end{equation}
But by Proposition \ref{prop:qivolume} and (\ref{eq:tilde:B:H:cap:U:star}),
\begin{equation}
\label{eq:different:slab:volume:comparasion} \omega_1^{-1}
\Vol'(Sl_2^1(H)) \le \Vol'(\phi(Sl_{2}^{1}(H)
\cap U_*)) \le \omega_1 \Vol'(Sl_2^1(H)),
\end{equation}
where $\omega_1$ depends only on $(\kappa,C)$. Now
(\ref{eq:measure:lower:bilip}) and
(\ref{eq:prelim:upper:measure:lip})  follow from
(\ref{eq:different:slab:volume:comparasion}),
(\ref{eq:slab:volume:comparison}) (\ref{eq:q:part15:bilishitz})
and the fact that the volume of a sufficiently coarsened
generalized slab is the area of the cross section times the
difference in height. \qed

\subsection{Proof of Lemmas \ref{lemma:measure:lower:bilip} and \ref{lemma:measure:upper:bilip}}

We first prove a preliminary estimate:

\begin{lemma}
\label{lemma:relation:h:plus:q} Given $p_1,p_2$ in $U^*$, then
$h^+(\phi(p_1),\phi(p_2))=q(h^+(p_1,p_2)) + O(\epsilon R)$.
\end{lemma}

\bold{Proof.}
By the definition of $U^*$ we can find $\tilde p_i$ in $U^*$ with
$\pi_{xz}(\tilde p_i)=\pi_{xz}(p_i)$ and vertical geodesic
segments $\gamma_i \subset U$ going up from $\tilde p_i$ which
come within $O(1)$ at $h^+(p_1,p_2)$.  Since each $\gamma_i$ is in
$U$, each $\phi(\gamma_i)$ is within $O(\epsilon R)$ of a vertical
geodesic $\tilde \gamma_i$ and the $\tilde \gamma_i$ come within
$O(\epsilon R)$ of one another only at
$h^+(\phi(p_1),\phi(p_2))+O(\epsilon R)$.  But by the definition
of the product map, $\tilde \gamma_1$ is within $O(\epsilon R)$ of
$\tilde \gamma_2$ at $q(h^+(p_1,p_2))$.
\qed\medskip

\bold{Proof of Lemma \ref{lemma:measure:lower:bilip}} Let $c_2 =
c_2(\epsilon,\theta)$ be a constant to be chosen later. Fix $i <
j$. Let $E_1 \subset S_X$ be such that for $x \in E_1$ there
exists a horocycle $H_x$ such that $x \in I_x \equiv \pi_-(H_x)$
and (\ref{eq:tilde:B:H:cap:U:star}) fails for some slab
$Sl_2^1(H_k)$ as in the statement of the lemma. Note that by
assumption $Sl_2^1(H_k) \subset B(R)$. Thus we have a cover of
$E_1$ by the intervals $I_x$. Then, by the Vitali covering lemma
there are intervals $I_k = \pi_-(H_k)$ such that the inequality
opposite to (\ref{eq:tilde:B:H:cap:U:star}) holds for $H_k$,
$\sum_k |I_k| \ge (1/5) |E_1|$, and also the $I_k$ are strongly
disjoint, i.e. for $j \ne k$, $d(I_j, I_k) \ge (1/2)
\max(|I_j|,|I_k|)$. Then the sets $SL_2^1(H_k)$ are also disjoint.
By construction, $|SL^1_2(H_k) \cap U_*^c| \ge c_2 |Sl^1_2(H_k)|$.
Summing this over $k$, we get that
\begin{displaymath}
|B(R) \cap U_*^c| \ge c_2 \sum_k |Sl_2^1(H_k)| \ge
(c_2/2) \sum_k
|h_1(H_k) - h_2(H_k)| |I_k| |S_Y|.
\end{displaymath}
Since $|B(R) \cap U_*^c| \le \theta R |S_X||S_Y|$, we get
\begin{displaymath}
|E_1| \le 5 \sum_k |I_k| \le \frac{10 \theta}{\beta' c_2} |S_X|.
\end{displaymath}

If $c_1c_2\beta'=20\theta$ this implies that
$|E_1|<\frac{c_1}{2}|S_X|$. So letting $E_{**}=S_X \backslash
E_1$, we are done.
\qed\medskip

\bold{Proof of Lemma \ref{lemma:measure:upper:bilip}} We construct
$E_*$ as a subset of $E_{**}$ from Lemma
\ref{lemma:measure:lower:bilip}, so any $H$ satisfying the
hypotheses of Lemma \ref{lemma:measure:upper:bilip} satisfies the
conclusions of Lemma \ref{lemma:measure:lower:bilip}.

 We now
show that $\tilde{Sl}_2^1(H)\subset\phi(B(R))$.  Recall that $H_0$
is more than $4 \kappa^2 \beta'' R$ from the edge of $B(R)$.  By
definition
$$\tilde Sl_{2}^{1}(H)=S( \cC_{q(h_c)}(f(\pi_-(H)), \cC_{g(h_2)}(g(S_Y{\cap} S)),
q(h_2), q(h_1)).$$

\noindent Since $S \subset U_*$, for any $y{\in}S_Y{\cap}S$ there
is a point $p=(x,y,z) {\in} S$ such that $\phi$ maps $p$ to within
$O(\epsilon R)$ of $(f(x),g(y),q(z))$ with $x$  in $\pi_-(H)$. The
point $p$ is at most $\beta'' R$ from $H_0$ and $\phi(p)$ is
within $O(\epsilon R)$ of a vertical geodesic $\gamma$ which stays
within $O(\epsilon R)$ of the image of a vertical geodesic through
$p$. Note that any point $q$ in $S(f(\pi_-(H)),g(S_Y \cap
S),q(h_2), q(h_1))$ is on a vertical geodesic $\gamma'$ which
stays within $\epsilon R$ of the image of a geodesic which passes
through $S$ and therefore through $H_0$.  The point $\phi(p)$ is
within $\kappa \beta'' R$ of where the geodesics $\gamma$ and
$\gamma'$ come within $O(\epsilon R)$ since $p$ is within $\beta''
R$ of the point where the corresponding geodesics in the domain
come close. This implies that $q$ is within $3.1 \kappa\beta''R$
of $\phi(p)$. By the definition of coarsening, this implies that
any point in $S( \cC_{q(h_c)}(f(\pi_-(H)),
\cC_{g(h_2)}(g(S_Y{\cap} S)), q(h_2), q(h_1))$ is within
$4\kappa\beta''R$ of $\phi(p)$. By our assumptions on $S$ and $p$,
this shows that $\tilde{Sl}_2^1(H)\subset\phi(B(R))$.

Let $c_3 = c_3(\epsilon,\theta,\beta')$ be a constant to be chosen
later. Let $E_2 \subset S_X\setminus E_1$ be such that for $x \in
E_2$ there exists a horocycle $H_x$ such that $x \in I_x \equiv
\pi_-(H_x)$ and (\ref{eq:hat:B:H:cap:N}) fails. Thus we have a cover
of $E_1$ by the intervals $I_x$. Then, by the Vitali covering lemma
there are intervals $I_k = \pi_-(H_k)$, such that the inequality
opposite to (\ref{eq:hat:B:H:cap:N}) holds for $H_k$ instead of $H$,
$\sum_k |I_k| \ge (1/5) |E_2|$, and also the $I_k$ are strongly
disjoint, i.e. for $l \ne k$, $d(I_l, I_k) \ge (1/2)
\max(|I_l|,|I_k|)$.

We now claim that
\begin{equation}
\label{eq:phi:Sh:H:k} \phi(\Sh(H_k,O(1))^c \cap U_*) \cap
\tilde{Sl}_2^1(H_k) = \emptyset.
\end{equation}
Indeed suppose $p \in \Sh(H_k,O(1))^c \cap U_*$, and $\phi(p) \in
\tilde{Sl}_2^1(H_k)$. By the definition of $\hat{Sl}_{2}^{1}(H_k)$,
$\pi_-(\phi(p)) \subset \cC_{q(h_c(H))}(f(\pi_-(H_k)))$. Hence, by
the definition of coarsening, there exists $p' \in \Sh(H_k,O(1))
\cap U_*$ such that $h^+(\phi(p),\phi(p')) = q(h_c(H)) + O(1)$.
Since $p_1 \in \Sh(H_k,O(1))^c$ and $p_2 \in \Sh(H_k,O(1))$, we have
$h^+(p_1,p_2)>h(H_k)+O(1)$. This contradicts Lemma
\ref{lemma:relation:h:plus:q}, and thus (\ref{eq:phi:Sh:H:k}) holds.
The same argument shows that the sets $\tilde{Sl}_2^1(H_k)$ are
disjoint.

Suppose $p \in \tilde{Sl}_{2}^{1}(H_k)$, $q(h_2(H_k))+O(\epsilon R)
< h(p) < q(h_1(H_k)) - O(\epsilon R)$, and $p \not\in N_{O(\epsilon
  R)}(\phi(Sl_2^1(H_k) \cap B(R)))$. We claim that $p \not\in
\phi(U_*)$. Indeed, if $p = \phi(p')$ where $p' \in U_*$, then by
(\ref{eq:phi:Sh:H:k}), $p' \not\in \Sh(H_k,O(1))^c$. But since
$h_2(H_k) < h(p') < h_1(H_k)$, we have $p' \in Sl_2^1(H_k) \cap
B(R)$. This is a contradiction, and hence $p \not\in \phi(U_*)$.
This implies that $\phi(U^* \cap Sl_2^1(H_k)^c) \cap
\hat{Sl}_2^1(H_k)$ contributes negligibly to the measure of
$\hat{Sl}_2^1(H_k)$, i.e. the contribution goes to zero as
$\epsilon$ goes to zero.  So to complete the proof, we need only
control $\Vol'(\tilde{Sl}_{2}^{1}(H_k) \cap \phi(B(R) \cap U_*^c)).$

Thus, since we are assuming the opposite inequality to
(\ref{eq:hat:B:H:cap:N}), we have $\Vol'(\tilde{Sl}_{2}^{1}(H_k)
\cap \phi(B(R) \cap U_*^c) \ge c_3 |\tilde{Sl}_{2}^{1}(H_k)|$. But
then, using the disjointness of the $\tilde{Sl}_{2}^{1}(H_k)$ we get
\begin{multline}
\label{eq:estimate:size:image:bad} \Vol'(\phi(B(R) \cap
U_*^c)) \ge c_3 \sum_k |\hat{Sl}_{2}^{1}(H_k)| \ge (c_3)(1-c)
\Vol'(\phi(Sl_2^1(H_k) \cap U_*)) \ge \\
\ge \omega_3 c_3 \sum_k |Sl_2^1(H_k) \cap U_*| \ge
\omega_4 c_3 \beta R \sum_k |I_k| |S_Y|.
\end{multline}
\noindent The first inequality is our assumption, the second uses
equation $(\ref{eq:slab:volume:comparison})$.  The third is
Proposition \ref{prop:qivolume} and also uses the fact that each
$I_k$ contains a point of $S_X\setminus E_1$ to conclude that
$|Sl_2^1(H_k) \cap U_*| \ge (1/2) |Sl_2^1(H_k)|$.

Since by Proposition \ref{prop:qivolume}, $\Vol'(\phi(B(R)
\cap U_*^c)) \le \omega_5 \theta R |S_X| |S_Y|$, we get
\begin{displaymath}
|E_2| \le 5 \sum_k |I_k| \le \frac{\omega_6 \theta }{c_3\beta}
|S_X|.
\end{displaymath}

And so $|E_2|<\frac{c_1}{2}$, provided $c_3 c_1 \beta =2\omega_6
\theta$.
So after letting $E_*=S_X \backslash E_1{\cup}E_2$, the proof
is complete.
\qed\medskip

\subsection{The map on heights}
\label{subsection:map:on:heights}

%Throughout this subsection, for simplicity of notation, we

Suppose $B(R) \subset X(n)$ is a box, and
suppose $\phi: B(R) \to X(n')$ is an $(\epsilon,\theta)$ almost-product
map. Then by definition, there exists a partially defined product map
$\hat{\phi} = (f,g,q)$ and a subset $U \subset B(R)$ with $|U| \ge
(1-\theta) |B(R)|$ such that
\begin{equation}
\label{eq:phi:close:hat:phi}
d(\phi|_U, \hat{\phi}) = O(\epsilon R)
\end{equation}

\begin{proposition}[Map on heights]
\label{prop:map:on:heights} Let $\beta \ll \beta' \ll \beta'' \ll
1$ be as in \S\ref{subsection:measure:bilipschitz}. Write $B(R) =
S_X \cross S_Y \cross [h_{bot}, h_{top}]$. Suppose the
following hold:
\begin{itemize}
\item $h_{bot} < z_{bot} < z_{top} < h_{top}$.
\item $4 \beta R \le |z_{top} - z_{bot}| \le \beta' R$.
\item $|h_{top} - z_{top}| > 4 \kappa^2 \beta'' R$.
\item $|z_{bot} - h_{bot}| > 4 \kappa^2 \beta'' R$.
\end{itemize}

Then
there exists a set $S{\subset}B(R)$ as in Lemma
\ref{lemma:measure:upper:bilip} and a function $\epsilon' =
\epsilon'(\epsilon,\theta)$ with $\epsilon' \to 0$ as $\epsilon
\to 0$ and $\theta \to 0$ such that for all $z \in
[z_{bot},z_{top}]$,
\begin{displaymath}
\label{eqn:heightaffine} q(z) = A z - \frac{1}{B_X'}\log
\frac{|\cC_{q(z_{bot})}(g(S_Y {\cap}S))|}{|S_Y|} + O(\epsilon' R),
\end{displaymath}
where $A = B_{X(n)}/B_{X(n')} = B_X/B_X'$ is the ratio of branching constants. In
particular, if $n = n'$, $A=1$.
\end{proposition}

\bold{Remark.}In all applications of Proposition
\ref{prop:map:on:heights}, we change $q$ by $O(\epsilon' R)$ in
order to have (\ref{eqn:heightaffine}) hold with no error term.

\bold{Remark.} For any $n$, $n'$ there exists a standard map
$\hat{\phi} = (f,g,q):
X(n) \to X(n')$ with $q(z) = A z$. For solvable
groups $\hat{\phi}$ is simply a homothety, for Diestal-Leader graphs
it is given by collapsing levels.
\bigskip

The rest of this subsection will consist of the proof of
Proposition~\ref{prop:map:on:heights}.  Apply
Lemmas~\ref{lemma:measure:lower:bilip} and
\ref{lemma:measure:upper:bilip} to get a set $E_* \subset S_X$.
Let $H,H_0$ be $y$ horocycles that satisfy the conditions of
Lemmas \ref{lemma:measure:lower:bilip} and
\ref{lemma:measure:upper:bilip} with $h(H_0)>h(H)>z_{top}$. In
particular $\pi_-(H)$ contains a point of $E_*$.  Choose an
arbitrary $z \in [z_{bot},z_{top}]$, let  $h_1 = z$, $h_2 =
z_{bot}$.  For the remainder of this subsection we simplify
notation by writing $g(S_Y)$ for $g(S_Y{\cap}S)$.

\begin{lemma}
\label{lemma:map:on:verticals} There exists a function $\epsilon'
= \epsilon'(\epsilon,\theta)$ with $\epsilon' \to 0$ as $\epsilon
\to 0$ and $\theta \to 0$ such that the following holds: Let $\cF$
denote the set of vertical geodesic segments in
$Sl_2^1(H)$, and let $\hat{\cF}$ denote the set of
vertical geodesic segments in $\tilde{Sl}_2^1(H)$. Then exists a
subset $\cF' \subset \cF$ with $|\cF'| \ge (1/2) |\cF|$ and a map
$\psi: \cF' \to \hat{\cF}$ which is at most $e^{\epsilon'
  R}$ to one. Also there exists a subset $\hat{\cF}' \subset
\hat{\cF}$ with $|\hat{\cF}'| \ge (1/2) |\hat{\cF }|$
and a map $\hat{\psi}: \hat{\cF}' \to \cF$ with is at most
$e^{\epsilon' R}$ to one.
Hence,
\begin{equation}
\label{eq:log:sizes:close}
\log |\cF| = \log |\hat{\cF}| + O(\epsilon' R).
\end{equation}
\end{lemma}

\bold{Proof.} Let $c_2$ be as in (\ref{eq:tilde:B:H:cap:U:star}).
We let $\cF'$ to be the set of vertical geodesics in
$Sl_2^1(H)$ more than $O(\epsilon R)$ from the edges and which spend
at least $1-\sqrt{c_2}$ fraction of their length in $U_*$. Then,
by (\ref{eq:tilde:B:H:cap:U:star}), $|\cF'| \ge (1/2) |\cF|$.
Now since $\phi$ is an almost-product map, for each $\gamma \in \cF'$
there exists a geodesic $\hat{\gamma} \in \hat{\cF}$ such that
$\phi(\gamma \cap U_*)$ is within $O(\epsilon R)$ of $\hat{\gamma}$.
We define $\psi(\gamma) = \hat{\gamma}$. The map $\psi$ is
at most $e^{O(\epsilon R + \sqrt{c_2}R)}$ to one since
two geodesics with the same image must be within $\epsilon R$ of
each other whenever they are in $U_*$, and by assumption there exist
points in $U_*$ on each geodesic within $O(\sqrt{c_2} R)$ of $h_{top}$
and $h_{bot}$.

The construction of the ``inverse'' map $\hat{\psi}$ is virtually
identical, except that one uses (\ref{eq:hat:B:H:cap:N}) instead of
(\ref{eq:tilde:B:H:cap:U:star}) and $c_3$ instead of $c_2$. In the
end, we can choose $\epsilon' = O(\epsilon + \sqrt{c_2} + \sqrt{c_3})$.
\qed\medskip

\begin{lemma}
\label{lemma:dilation:part}
For all $z \in [z_{bot},z_{top}]$,
\begin{displaymath}
q(z) - q(z_{bot}) = A(z - z_{bot}) + O(\epsilon' R)
\end{displaymath}
\end{lemma}

\bold{Proof.} We count vertical geodesics using Lemma
\ref{lemma:map:on:verticals}. Note that $$|\cF| \sim
|\pi_-(H)||S_Y|e^{B_X(h_1-h_2)},$$ and by
Lemma~\ref{lemma:generalized:slabs}, $|\hat{\cF}|$ is comparable
to
$$|\cC_{h_1}(f(\pi_-(H))||\cC_{h_2}(g(S_Y))|e^{B_{X'}(q(h_1)-q(h_2))},$$
where as above $h_1 = z$, $h_2 = z_{bot}$.
Then, by (\ref{eq:log:sizes:close}),
$$q(h_1)-q(h_2)= A(h_1-h_2) + \frac{1}{B_X'} \log
\frac{|\cC_{q(h_1)}(f(\pi_-(H)))||\cC_{q(h_2)}(g(S_Y))|}{|\pi_-(H)||S_Y|}
+ O(\epsilon' R).$$ Now by Corollary~\ref{cor:measure:lower:bilp}
the logarithm is bounded between two constants which depend only
on $\kappa$ and $C$. \qed\medskip

\bold{Proof of Proposition~\ref{prop:map:on:heights}.} Choose $h_1
= (z_{top}+ z_{bot})/2$, $h_2 = z_{bot}$. By
Lemma~\ref{lemma:measure:lower:bilip} there exists a horocycle
$H'$ with $h(H') = z_{top}$ so that (\ref{eq:measure:lower:bilip})
and (\ref{eq:prelim:upper:measure:lip}) hold for $H'$. Then
\begin{equation}
\label{eq:log:ratio:temp} \log
\frac{|\cC_{q(h_1)}(f(\pi_-(H')))||\cC_{q(z_{bot})}(g(S_Y))|}{|\pi_-(H')||S_Y|}
= O(1).
\end{equation}
By Lemma~\ref{lemma:relation:h:plus:q}, equation
(\ref{eq:measure:shadow:point}) and the fact that we coarsen below
the horocycle, we see that,
\begin{displaymath}
\log \frac{|\cC_{q(h_1)}(f(\pi_-(H')))| e^{-B_X' q(h(H'))}}{|\pi_-(H')|
  e^{-B_X h(H')}} = O(\epsilon R).
\end{displaymath}
Since $h(H') = z_{top}$, after rearranging we get,
\begin{displaymath}
\frac{1}{B_X'}\log \frac{|\cC_{q(h_1)}(f(\pi_-(H')))|}{|\pi_-(H')|} =
q(z_{top}) - A z_{top} + O(\epsilon R).
\end{displaymath}
Substituting into (\ref{eq:log:ratio:temp}) we get,
\begin{equation}
\label{eq:almost:done:here} \frac{1}{B_X'}\log
\frac{|\cC_{q(z_{bot})}(g(S_Y))|}{|S_Y|} = A z_{top} -
q(z_{top}) = A z_{bot} - q(z_{bot}),
\end{equation}
where we have used Lemma~\ref{lemma:dilation:part} for the last equality.
Now Proposition~\ref{prop:map:on:heights} follows from
(\ref{eq:almost:done:here}) and Lemma~\ref{lemma:dilation:part}.
\qed\medskip

\section{Proof of Theorem \ref{theorem:main:step:two}}
\label{section:thework}

In this section we prove Theorem \ref{theorem:main:step:two}. The
basic strategy is to show that for most horocycles $H$
intersecting $\phi{\inv}(B(L'))$, the image $\phi(H)$ is within
$\epsilon R$ of a horocycle, at least for most of it's measure.
This argument occupies the first four subsections.  Subsection
\S\ref{sec:subsection:end:proof:theorem2.1} completes the proof in a
manner analogous to \cite[Section 5.4]{EFW1}.

A key ingredient in our proofs is Lemma
\ref{lemma:illegal:1:circuit}.  The reader should think of this
``illegal circuit lemma" as a generalization or strengthening of
the ``quadrilaterals lemma" \cite[Lemma 3.1]{EFW1}.  The greater
generality comes from making weaker assumptions on the paths
forming the ``legs" of the ``quadrilateral".  Lemma
\ref{lemma:illegal:1:circuit} is used much like \cite[Lemma
3.1]{EFW1} to show that points along a horocycle must map by
$\phi$ to points approximately along a horocycle.

\bold{Recommendation to the reader:} We recommend that the reader
read this section first assuming that, for each $i{\in}I_g$, the
map $\phi$ restricted to $U_i$ in $B_i(R)$ is within $O(\epsilon
R)$ of a $b$-standard map.  Under this assumption, the
construction of the $\hat{S}$-graph can be omitted since it
suffices to consider only the $S$-graph.  The reader will find
that proofs in \S\ref{subsection:averaging:H:graph} and
\S\ref{subsection:circuits} simplify somewhat under this
hypothesis, but the main arguments in \S\ref{section:families}
remain essentially the same.

The primary difficulty that occurs here in dropping the assumption
that $\phi_i|_{U_i}$ is within $O(\epsilon R)$ of a $b$-standard
map is in guaranteeing that the map preserves the ``divergence
conditions" on pairs of vertical geodesics required to control
paths by the methods of \S\ref{subsection:circuits}.

\subsection{Constructing the $\hat{S}$ graph and the $H$-graph}
\label{subsection:build:H:graph}

Given a ``good enough" horocycle $H$ mostly contained in
$\phi{\inv}(B(L'))$, in this section we construct a graph which we
use to control $\phi(H)$.  To begin, we choose constants and make
precise the notion of a ``good enough" horocycle.

\bold{Choosing Constants:}
Let $\phi: X(n) \to X(n')$ be a $(\kappa,C)$ quasi-isometry. Chose
$\rho_1 \GG C$, and discretize on scale $\rho_1$ as described in
subsection \ref{subsection:discretization}. Let $B_X$ (resp.
$B_{X'}$) be the branching constant of the resulting graph and let
$B=\max\{B_X,B_{X'}\}$. Let $\epsilon > 0$ and $\theta > 0$ be
constants to be specified below, and let
$L'$ be sufficiently large so that
Theorem~\ref{theorem:partIinImage} applies, and fix a box $B(L')$.
We call the graph that is the discretization of $B(L')$ the
$S$-graph.

We now apply Theorem~\ref{theorem:partIinImage} to $B(L')$. We fix
$\epsilon \ll \alpha \ll \beta \ll \beta' \ll \beta''$ and apply
the arguments described in Section \ref{sec:measure:bilipshitz} to
each box $B_i(R)$ for $i{\in}I_g$ as in the conclusion of Theorem
\ref{theorem:partIinImage}, to obtain sets $(E_*^+)_i \subset
\partial^+X$ and $(E_*^-)_i \subset
\partial_-X$.
After replacing the set $U_i$ from
Theorem~\ref{theorem:partIinImage} with a slightly smaller set, we
can make sure that for all $(x,y,z) \in U_i$, $x \in (E_*^+)_i$,
$y \in (E_*^-)_i$. We still have $|U_i| \ge (1-\delta_0)|B_i(R)|$,
where $\delta_0 \to 0$ as $\epsilon \to 0$ and $\theta \to 0$.  As
remarked following Proposition \ref{prop:map:on:heights}, we
further modify $q_i$ so that it satisfies (\ref{eqn:heightaffine})
with no error term. This makes $\hat{\phi_i}$ within $O(\epsilon'
R)$ of $\phi$ where $\epsilon'$ goes to zero as $\epsilon \to 0$
and $\theta \to 0$.

We then choose $0<\eta \ll 1$ such that $\rho_1 \ll 1/\eta$ (We
mean that for any function $f(\rho_1)$ and any quantity $u$ which
is labeled $O(\eta)$ in the argument, $f(\rho_1)$ is much less
than $1$.)

We then choose $\rho_2 \GG \rho_1$ so that $f(\rho_1)/B^{\rho_2}
\ll \eta$, where $f(\rho_1)$ is any function of $\rho_1$ which
arises during the proof. Now pick $\rho_3$, $\rho_4$, $\rho_5$ so that
$\rho_2 \ll \rho_3 \ll \rho_4 \ll \rho_5$.

Choose $0 < \delta_0 \ll 1$, so that $\rho_5 \ll 1/\delta_0$
(The last inequality means that for any function $f(\rho_5)$
and any function $g(\delta_0)$ going to $0$ as $\delta_0 \to 0$
which arise during the
argument, $f(\rho_5) g(\delta_0) \ll 1$. We also assume that
$1/\eta \ll 1/\delta_0$ (i.e. for any quantity $u$ labeled $O(\eta)$
and any function of $g(\delta_0)$ going to $0$ as $\delta_0 \to 0$
which arises during the proof, we have $f(\delta_0)  \ll u$.

%\mc{remark: $\delta_0$ is essentially the $\theta$ from
%  Theorem~\ref{theorem:partIinImage}. }

% We know that for any $\delta_0$ there exist $epsilon R < R < L$ such that if
% we take any box $B(L)$, and tile it by boxes of size $R$, for at least
% $(1-\delta_0)$ fraction of these boxes there exists a ``good''
% subset of the
% box of relative measure at least $(1-\delta_0)$, such that the
% restriction of $\phi$ to the subset is within $epsilon R$ of a standard map.

\bold{Recap.}
We have
\begin{displaymath}
C \ll \rho_1 \ll \rho_2 \ll \rho_3 \ll \rho_4 \ll \rho_5 \ll (1/\delta_0).
\end{displaymath}
also
\begin{displaymath}
\rho_1 \ll 1/\eta \ll 1/\delta_0
\end{displaymath}
and
\begin{displaymath}
\rho_5 \ll \epsilon' R \ll R \ll L'.
\end{displaymath}
We do not assume that e.g. $e^{\epsilon' R} \delta_0$ is small.

\bold{Note.} We assume $e^{\epsilon' R} \GG L'$. Both of these are
consequences of the proof of Theorem~\ref{theorem:partIinImage}.
We always assume that any path we consider has length $O(L')$
which is much smaller then $e^{\epsilon'
  R}$.
\medskip

\bold{The sets $U'$ and $U$.} Let $U_i$, $i \in I_g$ be as in the
second paragraph of this subsection.  Let $U' = \bigcup_i U_i$.
Then $|{U'}^c \cap \phi^{-1}(B(L'))| \le 2 \delta_0
|\phi^{-1}(B(L'))|$.

Let $U''$ denote the subset of $\phi^{-1}(B(L'))$ which is
distance at most $\rho_1$ from $U'$. Then $|U''| \ge |U'| \ge (1 -
\delta_0') |\phi^{-1}(B(L'))|$. Also note that since $\rho_1 \ll
\epsilon' R$, for $i \in I_g$, the restriction of $\phi$ to $U
\cap B_i(R)$ is an  $(\epsilon' R +\rho_1,\theta)$-almost product
map. We define a set $U$ by $U^c=N_{\rho_5+\rho_1}({U''}^c)$.  An
elementary covering lemma argument shows that $|U| \ge
(1-\delta_0'')|\phi^{-1}(B(L'))|$ where $\delta_0''$ goes to zero
with $\delta_0'$.

\bold{Favorable Horocycles.} We define a horocycle $H$ to be {\em
favorable} if $H$ does not stay within $\beta'' R$ of the walls of the
$B_i(R)$, and also
$$|H \cap U'| \ge (1 - \delta_0''') |H \cap
  \phi^{-1}(B(L')).$$
  \noindent We call $H$ {\em very favorable} if the same holds with
$U$ in place of $U'$.

\bold{Remark:} If a horocycle is very favorable, any
horocycle within $\rho_5$ of it is favorable.

\begin{lemma}
\label{lemma:most:favorable} There exists $\hat{\theta} > 0$ such that
the fraction of $B(L')$ which is contained in the image of a
very favorable $x$-horocycle and a very favorable $y$-horocycle is at least
$(1-\hat{\theta})$. Here $\hat{\theta}$ is a function of $\delta_0$
and $\beta''$  which goes to $0$ as $\delta_0 \to 0$ and $\beta'' \to 0$.
\end{lemma}

\bold{Proof.} This is immediate from the construction.
\qed\medskip

\bold{Notation.} For most of the argument, we fix
a very favorable horocycle $H$, whose image $\phi(H)$ intersects $B(L')$.
For notational simplicity, we assume that $H$ is an $y$-horocycle.
We also fix a favorable horocycle $H_0$ so that $\rho_5/2 < d(H,H_0) <
\rho_5$ and $H \subset \Sh(H_0, \rho_1)$.

\bold{The sets $I_g(H)$, $\tilde{B}$ and $U_*$.} Let $I_g(H)$
denote the set of indices $i \in I_g$ such that
\begin{equation}
\label{eq:good:box}
|H \cap U' \cap B_i(R)| \ge (1 - \delta_0'') |H \cap B_i(R)| > 0.
\end{equation}
Now let $\tilde{B} = \bigcup_{i \in I_g(H)} B_i(R)$, and let
$U_* = U' \cap \tilde{B}$.

\bold{Good and bad boxes.}
We refer to boxes $B_i(R)$ with $i \in I_g(H)$ as ``good boxes'', and
to boxes $B_i(R)$ intersecting $H$ with $i \in I \setminus I_g(H)$
as ``bad boxes''.

% We fix an favorable horocycle $H$ and assume in what follows that
% $H$ is an $x$-horocycle.  The argument for $H$ a $y$-horocycle is
% identical up to ``flips".

\bold{Shadows of $H$ and $\phi(H)$:} Let $h_1=h(H)-(\alpha+\beta)
R$ and $h_2=h(H)-(\alpha+ \beta +\frac{\beta'}{2})R$.
%In each box $B_i(R)$ intersecting $H$, we can consider the set
%$Sl_2^1(H)_i$. For any good $B_i(R)$, we can also consider
%$\hat{Sl}_2^1(H)_i$. Furthermore
For each $i{\in}I_g$, we let $h^i_0$ to be specified below be such
that $(\alpha+\frac{\beta}{2})R<|h(H)-h^i_0|<(\alpha+{\beta})R$.
For each $B(R)_i$ intersecting $H$, we denote
${W}(H)_i=h{\inv}(h_0){\cap}\Sh(H_0,\rho_1)$. For all bad boxes,
we fix $h_0^i=h(H)-(\alpha+\beta)R$, for good boxes $h_0^i$ will be
fixed during the proof of Lemma \ref{lemma:shadow:vertices} below.
For each $B_i(R)$ intersecting $H$ with $i{\in}I_g$ we let
$$\hat{W}(H)_i=\{(x,y,z)| x{\in}\cC_{q(h_1)}(f(\pi_-(H_0)), y{\in}\cC_{q(h_2)}(g(S_Y{\cap}Y), z=q(h_0)\}.$$

\noindent Let $W(H)=\cup_{i{\in}I}{W}(H)_i$, and $\hat{W}(H) =
\cup_{i{\in}I_g}\hat{W}(H)_i$. We define these sets in terms of
$H_0$ not $H$ so as to be able to consider points above $H$ in
certain arguments below. Recall $q$ is fixed so that the
$O(\epsilon' R)$ term in Proposition~\ref{prop:map:on:heights} is
$0$.  We let $R_i'=h(H)-h^i_0$.  We frequently suppress reference
to $i$ in our notation for $R'$ and $h_0$. 
%{\sc FIX $H$ VS $H_0$ below.}

\bold{Shadow vertices.} We now define a set of {\em shadow
vertices} in the discretization of $X(n')$.  By shifting the
discretization, we can assume that $\hat{W}(H)$ contains a
$\rho_1$ net of $S$-vertices.
 Every $S$-vertex in $\hat{W}(H)$ is a
{\em shadow vertex}. If some vertical geodesic going down $\beta'
R$ from $s$ contains a point of $\phi(U')$ below $h_1$ and $s$ is
not within $10 \kappa \epsilon' R$ of an edge of $\hat{W}(H)$ then
we call $s$ a {\em good shadow vertex}. Any $S$-vertex in
$\hat{W}(H)$ which is not a good shadow vertex is a {\em bad
shadow vertex}. We now add additional shadow vertices, not
necessarily in $\hat{W}(H)$. We also make any $S$ vertex in
$N_{\rho_1}\phi({U_*}^c{\cap}W(H))$ a {\em bad shadow vertex},
even if it is a good shadow vertex by our previous definition. The
bad shadow vertices in $N_{\rho_1}\phi({U_*}^c{\cap}W(H))$ are not
necessarily close to $\hat{W}(H)$, even if they come from good
boxes.  While these bad shadow vertices are not well controlled,
they make up a small proportion of all shadow vertices and so do
not interfere with our arguments, see
Lemma~\ref{lemma:shadow:vertices} below.

For either good boxes or bad boxes, the number of shadow vertices
coming from $B_i$ is proportional to the length of $H{\cap}B_i$.
The proportionality constant depends only on $\kappa,C$ and the
geometry of the model spaces.

\begin{lemma}
\label{lemma:shadow:vertices}   There is a constant
$c_4(\delta_0,\epsilon')$ such that, for appropriate choices of
$h_0^i$,  only $c_4$ fraction of all shadow vertices are bad and
$c_4$ goes to zero as $\delta_0,\epsilon'$ go to zero.
\end{lemma}

\bold{Proof.} Bad shadow vertices are defined in two stages.  First
we have the set $S_1$ of vertices in $\hat{W}(H)_i$ not within
$10\kappa\epsilon' R$ of an edge and not within $\beta R$ of a
point in $\phi(U')$ below $h_0$.  That this set has small measure
in $\hat{W}(H)_i$ follows from two facts. First, the subset within
$10\kappa\epsilon' R$ of the boundary is has measure going to zero
with $\epsilon'$. Second, if $S_1$ contains $\theta$ fraction of
the vertices in $\hat{W}(H)_i$, then the set of points on
geodesics going down $\beta' R$ from $S_1$ contains
$\frac{\theta}{2}$ fraction of the measure of
$\tilde{Sl}_2^1(H_0)$. But $\phi(U')$ contains $1-c_4$ of the
measure in $\tilde{Sl}_2^1(H_0)$ by Lemma
\ref{lemma:measure:lower:bilip}, so this implies that
$\theta<2c_4$.  Since $c_4$ goes to zero with $\epsilon'$ and
$\delta_0$, this implies $|S_1|$ goes to zero with them as well.

In the second stage, we enlarge the set of bad vertices in
$\hat{W}(H)$ by adding the set $N_{\rho_1}\phi(U'^c{\cap}W(H)))$
to the set of bad vertices.  The fraction of shadow vertices
coming from bad boxes, being proportional the the fraction of the
length of $H$ in bad boxes, clearly goes to zero as $\delta_0$
goes to zero.  So it suffices to control the size of
$N_{\rho_1}\phi(U_*^c{\cap}W(H)_i))$. To show that
$W(H)_i{\cap}U_*^c$ contains a small fraction of the measure of
$W(H)_i$, we use the flexibility in our choice of $h^i_0$. Here we
make this flexibility explicit by letting $W(H)_i(h^i_0)$ be the
set of possible $W(H)_i$'s, parametrized by choices of $h^i_0$.
Let $\rho(h^i_0)$ be the fraction of  $W(H)_i(h^i_0)$ contained in
$W(H)_i{\cap}U_*^c$. Consider the slab $Sl_2^1(H_0)$ with
$h_1=(\alpha+\frac{\beta}{2})R$ and $h_2=h(H)-(\alpha+ 2\beta)R$.
Since all $W(H)_i(h^i_0)$ are contained in $Sl_2^1(H)$, using
Lemma \ref{lemma:measure:lower:bilip}, we have that
$$\sum_{h_0=(\alpha+\frac{\beta}{2})R}^{(\alpha+\beta)R}\rho(h_0)\leq
2c_3(\epsilon',\theta)$$ \noindent which implies that for some
$h^i_0$, we have $\rho(h_0)<2\sqrt{c_3}$.  We fix some $h_0^i$
with this property.

Lastly we need to see that this contribution remains small
relative to the number of good shadow vertices coming from
$B_i(R)$.  To see this, we use Corollary
\ref{cor:measure:lower:bilp} which implies that $|W(H)_i| \sim
|\hat{W}(H)_i|$ for constants depending only on $\kappa$ and $C$.
Combined with Proposition \ref{prop:qivolume}, this implies that the
ratio of $|N_{\rho_1}\phi(U_*^c{\cap}W(H)_i))|$ to the number of
good vertices in $\hat{W}(H)_i$ goes to zero with $\epsilon'$ and
$\delta_0$.

\qed

\bold{The $\hat{S}$-graph.} It is convenient to modify the
$S$-graph near the image of $H$. For $x \in \partial_-X$, $y \in
\partial^+X$ and $t \in [q(h_0),q(h(H_0))-\frac{\rho_5}{4}]$, let
$\gamma_{x,y}(t) = (x,y,t)$, so that $\gamma_{x,y}$ is a vertical
geodesic segment of length $q(h(H_0))-\frac{\rho_5}{4}-q(h_0)$.
Let $K_i$ be the union of $\gamma_{x,y}$ where
${x,y,q(h_0)}{\in}\hat{W}(H)_i$.  We begin by replacing $K_i$ as a
subset of the $S$ graph by the disjoint union of the
$\gamma_{x,y}$. We then define the $\hat{S}$ graph by defining a
new set of vertices and a new incidence relation on $K_i$. For $1
\le j \le q(h(H_0))-\frac{\rho_5}{4}-q(h_0)/\rho_1$ let $t_j =
q_i(h_0) + j \rho_1$.  We introduce {\em pre-vertices} along each
$\gamma_{x,y}$ at each $t_j$.  An {\em irregular} $\hat{S}$-vertex
will be an equivalence class of pre-vertices. Each pre-vertex has
coordinates $\{x,y,t_j\}$ At each height level $t_j$ in $X'(n)$,
we tile the $y$-horocycle by by disjoint segments $T_y$ of length
$10\rho_1$. At each height level $q{\inv}(t_j)$ in $X(n)$ we tile
each $x$ horocycle by disjoint segments $T_x$ of length
$10\kappa^2\rho$. (These tilings are best thought of as tilings of
horocycles in the corresponding trees or hyperbolic planes.)  We
identify two pre-vertices if:

\begin{enumerate}
\item their projections to the $yt$ plane are in the same $T_y$
and \item the points $(f_i{\inv}(x),q_i{\inv}(t_j))$ and
$(f_i{\inv}(x'),q_i{\inv}(t'_j))$ are in the same $T_x$. \item
$\pi_-(T_x) \cap f((E_**^-)_i)$ contains at least half the measure
in $\pi_-(T_x)$.
\end{enumerate}

Any segment ending at a bad shadow vertex is removed.  The
$\hat{S}$-vertices which are $S$-vertices outside of $K_i$ are
called {\em regular}.

% Let $K_i \subset X(n')$ denote the union of all the
%vertical segments $\gamma_{x,y}$ of length $q_i(h(H_0)) - q_i(h_3)$
%going up passing through $\hat{W}(H)_i$.  The $\hat{S}$-graph will be
%identical to the $S$-graph outside the union of the $K_i$. The
%$\hat{S}$-graph inside $K_i$ is obtained by taking the segments
%$\gamma_{x,y}$ and for each $j$
%attaching together to form a vertex those that are incident at height $t_j$
%(as defined in the previous paragraph). The $\hat{S}$-vertices
%obtained this way are called {\em irregular}; the $\hat{S}$-vertices
%which are also $S$-vertices are called {\em regular}.

\bold{The cloud of an $\hat{S}$-vertex.} Note that for any
$\hat{S}$-vertex $v$, $h(v)$ and the $y$-coordinate of $v$ are well
defined. For an irregular $\hat{S}$-vertex the $x$ coordinate is
``fuzzy''. More precisely, the {\em cloud} of an $\hat{S}$-vertex
$v$ is the set of points at height $h(v)$ which are on the vertical
segments incident to $v$. Then for a regular $\hat{S}$-vertex, the
cloud is essentially a point (it has size $O(\rho_1)$), whereas for
an irregular $\hat{S}$-vertex the cloud can have size $D\epsilon' R$
where $D$ is a constant depending only on $\kappa, C$ and the model
geometry.

\bold{The set $\hat{\phi}(H')$.}
Note that if $H'$ is within $\rho_4$ of $H$, then the set
$\hat{\phi}_i(H')$ consisting of the $\hat{S}$-vertices $v$ with
$q_i(h(H')) = h(v)+O(\rho_1)$ and the $x$-coordinate of $H'$ is
$f_i^{-1}(v)+O(\rho_1)$ is well defined. (The notation is explained by
the fact that for any $v \in \hat{\phi}_i(H')$, $\hat{\phi}_i^{-1}(v)$
is within $O(\rho_1)$ of $H'$.)
We then define $\hat{\phi}(H')$ to be $\bigcup_{i \in I_g} \hat{\phi}_i(H')$.

\begin{lemma}
\label{lemma:hat:S:constant:valence}
There exist constants $M_l$ and $M_u$ depending only on $\kappa, C$
such that for any two $\hat{S}$-vertices
$v_1$ and $v_2$ in $B(L')$, the ratio of the number of vertical geodesics in
$B(L')$ passing through $v_1$ to the number of vertical geodesics in
$B(L')$ passing through $v_2$ is bounded between $M_l$ and $M_u$.
\end{lemma}

\bold{Proof.} The proof is mainly a computation of the valence of
(i.e. the number of vertical paths incident to) an irregular
vertex. We give the proof in the $\DL$ case first. In the $\DL$
case the valence of a regular vertex is clearly $e^{{B_X'}L'}$ and
we will see that irregular vertices have the same valence.   For
$\Sol$, the valence of regular vertices can vary by a factor of
$2$ due to edge effects.  This same factor of $2$ occurs in the
first step of the computation below.

Let $h_{top}$ denote the height of the top of $B(L')$, and
$h_{bot}= h_{top} - L'$ denote the height of the bottom of
$B(L')$. Suppose $v$ is an irregular vertex in $K_i$. We can
choose a horocycle $H'$ so that $v \in \hat{\phi}_i(H')$, hence
$q_i(h(H')) = h(v)$. Note that by definition, $\pi_-(H')$ contains
a point in $E_{**}$. Then the number of paths going up from $v$ to
the height $h_{top}$ is $\approx e^{B_X'(h_{top} -
  h(v))}$.
Now the number of paths going down from $v$ to $\hat{W}(H)$ (at
height $h_0$) is
\begin{align*}
 \approx\, & |\cC_{q_i(h_0)}(f_i(\pi_-(H')))| e^{-B_X' q_i(h_0)} &\quad
 &\text{ by
  (\ref{eq:measure:shadow:point})} \\
 \approx\, & |\cC_{q_i(h_0)}(f_i(\pi_-(H')))| e^{-B_X h_0}
 \frac{|\cC_{q_i(h_2)}(g_i(S_Y \cap S))|}{|S_Y|}
 &\quad &\text{ by Proposition~\ref{prop:map:on:heights}} \\
 \approx\, & |\pi_-(H')| e^{-B_X h_0} &\quad &\text{ by
  Corollary~\ref{cor:measure:lower:bilp}} \\
 \approx\, & e^{B_X(h(H') -h_0)} & \quad &\text{ by
   (\ref{eq:measure:shadow:point})} \\
  =\, & e^{B_X'(h(v) - q_i(h_0))} & \quad &\text{ by
   Proposition~\ref{prop:map:on:heights}}
\end{align*}
Thus the total number of paths going down from $v$ to $h_{bot}$ is
$\approx e^{B_X'(h(v) - h_{bot})}$, and thus the total number
of paths incident to $v$ is $\approx e^{B_X' L}$ as required.

% {\sc FIX FOLLOWING ARGUMENT. GIVE MORE DETAIL.}
% We recall that $\hat{SL}_2^1(H)$ is a disjoint union of boxes. For
% any point in $h{\inv}(q(H)){\cap}\hat{SL}_2^1(H)$ every monotone
% geodesic going towards $\hat{W}(H)$ is contributes to the valence.
% In view of Lemma~\ref{lemma:measure:lower:bilip},
% Lemma~\ref{lemma:measure:upper:bilip} and
% Proposition~\ref{prop:map:on:heights}, there
% are between $M_l \exp({{B_X'}R''})$ and $M_u \exp({B_X'}R'')$ such
% geodesics leaving any such point.  {\sc EXPLAIN BETTER}
% Counting with
% appropriate multiplicity, it is easy to that this implies that the
% valence of a good $H$-vertex is between $M_l \exp({{B_X'}L'})$ and
% $M_u \exp({{B_X'}L'})$.

%For $\Sol$, there is an additional error introduced by edge
%effects.  Even for bad vertices, the valence is between
%$\frac{\exp({{B_X'}L'})}{2}$ and $\exp({{B_X'}L'})$, where the
%lower number is achieved exactly by vertices at the edge of
%$B(L')$.  A similar edge effect occurs when counting the number of
%vertical paths leaving a point of
%$h{\inv}(q(H)){\cap}\hat{SL}_2^1(H)$.  The point is that
%$\hat{SL}_2^1(H)$ is a union of boxes of fixed size and points
%near the edge of those boxes may have half as many vertical
%geodesics leaving them in $\hat{SL}_2^1(H)$.  (This actually
%occurs only near the actual edges of $\hat{SL}_2^1(H)$). In any
%case, for $\Sol$, the valence is between $\frac{M_l}{2}
%\exp({{B_X'}L'})$ and $M_u \exp({{B_X'}L'})$.
\qed\medskip

\bold{The $H$-graph.}
An irregular $\hat{S}$-vertex $v \in K_i$ is an $H$-vertex if and only if
$q_i(h(H)) = h(v)+O(\rho_1)$ and the $x$-coordinate of $H$ is
$f_i^{-1}(v)+O(\rho_1)$. These vertices are called ``good''.
(Note that for any good $H$-vertex $v \in K_i$, $\hat{\phi}_i^{-1}(v)$
is within $O(\rho_1)$ of $H$).

We also declare the ``bad'' $H$-vertices to be the bad shadow vertices
(these are always regular $\hat{S}$-vertices. The ``good'' and ``bad''
vertices thus defined comprise all the vertices of the $H$-graph. An
edge of the $H$-graph is a vertical path in the $\hat{S}$-graph which
either connects two $H$-vertices, or connects an $H$-vertex to the top
or bottom of the box $B(L')$.  An edge with one endpoint at the top or
bottom of the box is called an {\em leaf edge}.

We will count edges with multiplicity.  An edge has multiplicity
equal to the number of vertical paths in the $\hat{S}$-graph which contain it.

\medskip

%\bold{Warning:} the distance from $\phi(p)$ to the corresponding vertex is
%$O(r)$ for a good vertex. For a bad vertex, it is $O(r)$ for one copy,
%and $O(\rho_1)$ for the other possible copy.
%\medskip

%\bold{Note.} If $i \in I_g(H)$ so that
%$B_i(R) \subset \tilde{B}$, then every $S$-vertex on
%$\hat{\phi}_i(H \cap B_i)$ is an $H$-vertex.

%\bold{Recap.}
%Since $\beta^{\rho_4} \sqrt{\delta_0'} \ll 1$, and with a
%suitably chosen $\delta_0''$,
%the number of
%bad vertices is much smaller then the number of good vertices, and the
%proportion goes to $0$ as $\delta_0 \to 0$.
%\medskip

%\bold{Subdividing/adding some edges:} We now subdivide any edge
%that passes through an $H$ vertex in it's interior.  This
%occurs exactly when an edge passes through a guardian vertex. We
%have also added some bad vertices when modifying the
%$H$-graph.  We add (and subdivide) edges which are
%monotone rays incident (genuinely, not quasi) on any new bad
%vertex added at that step.

%We define the valence of an $H$ vertex to be the number of
%edges in it's equivalence class.

\bold{Notation.} We denote the $H$-graph by $\cG(H)$. Let $\cV$
denote the set of vertices of $\cG(H)$, and let $\cE$ denote the
set of edges. Let $\cV_1 \subset \cV$ denote the set of ``good''
vertices as defined above.  We call an $H$ vertex $y$ oriented
(resp. $x$ oriented) if the horocycle segment containing it is a
$y$ horocycle (resp. $x$ horocycle).  We also refer to an
orientation for $\hat{S}$ vertices, which is just the orientation
of $H$ vertices in the same box.

\begin{lemma}
\label{lemma:uniformvalence} The valence of $H$ vertices is
bounded between two constants $M_l$ and $M_u$ depending only on
$\kappa,C$ and the model geometries. Furthermore $|\cV_1| \geq
(1-c_5)|\cV|$ where $c_5 = c_5(\epsilon',\delta_0)$ goes to zero
with $\delta_0 \to 0$ and $\epsilon' \to 0$.
\end{lemma}

\bold{Proof.} The first statement of the lemma is immediate from
Lemma~\ref{lemma:hat:S:constant:valence}. To show the final claim,
let $F$ denote the set of vertical paths passing through the good
shadow vetices. By definition, every such path is incident to a
good $H$-vertex, and also every vertical path incident to a good
$H$-vertex belongs to $F$. Thus $F$ is also equal to the set of
vertical paths incident to good shadow vertices. Let $A$ denote
the set of good shadow vertices. Since the valence of each
$H$-vertex is between $M_l$ and $M_u$ times the valence of each
good shadow vertex, we have $M_l |A| \le |F| \le M_u |A|$, and
$M_l |\cV_1| \le |F| \le M_u |\cV_1|$. Thus, $|\cV_1| \ge
(M_l/M_u)^2 |A|$. But by Lemma~\ref{lemma:shadow:vertices}, $|\cV
\setminus \cV_1| \le c_4 |A|$, where $c_4(\epsilon',\delta_0) \to
0$ as $\epsilon' \to 0$ and $\delta_0 \to 0$. Thus $|\cV \setminus
\cV_1| \le (M_l/M_u)^2 c_4 |\cV_1|$. \qed\medskip

%\bold{The sets $W_k$, $W_k'$, $\cW_k$  and $\cW_k'$.}
%For each integer $k$, $0 < k \le R/\rho_1$,
%let $W_k$ denote the set (in the domain)
%obtained by walking down for $k$ edges from $H$, and let $\cW_k$
%denote the image of $W_k$ under $\phi$. Note that for a suitable $k
%\approx R/\rho_1$, $W_k = W(H)$. Let
%$\cW_k$ denote the image of $W_k$ under $\phi$.

%Let $W_k' \subset W_k$ denote the set obtained by walking down for $k$
%edges from a point in $H \cap B_i(R)$ where $i \in I \setminus
%I_g(H)$, i.e. from a point on $H$ in one of the ``bad boxes''.
%Let $\cW_k'$ denote the image of $W_k$ under $\phi$.

%Note that $\cW_0'  \subset \cV
%\setminus \cV_1$, and the cardinality of $W_k'$ (and hence
%approximately that of $\cW_k'$) is independent of $k$.

\subsection{Averaging over the $H$-graph}
\label{subsection:averaging:H:graph}

In order to make our geometric arguments in the next section, we need to show
that the paths and configurations we consider in the $H$ graph do not involve any
bad shadow vertices. I.e. to show that these paths and configurations only come near
the horocycle in the good set, at places where we have control over the map $\phi$.  In this section we
use multiple applications of the Vitali covering lemma in order to guarantee that ``most"
configurations in the $H$-graph avoid bad shadow vertices.  A key fact is
that while neither $\Sol$ nor $\DL(m,m)$ satisfy the sort of doubling condition needed for the
Vitali covering lemma, the space of vertical geodesics does.

 Choose $0 < \theta_3 <
\theta_4 \ll 1$. The $\theta$'s will be functions of $\delta_0$
which go to $0$ as $\delta_0 \to 0$.

\begin{definition}[Good Edges]
\label{def:good:edges}
The following defines sets of ``good'' edges. See also
Definition~\ref{def:good:vertices}.
\begin{itemize}
\item[{$\cE_1$:} ] Either connects two vertices in $\cV_1$ or is a leaf
  edge based on a vertex of $\cV_1$.

\item[{$\cE_3$:} ] An $\cE_1$ edge $e$ such that for for all
  $\hat{S}$-vertices $x \in
  e$, $1-\theta_3$ fraction of the edges (forward) branching at $x$ are in
  $\cE_1$. (note that $x$ is not supposed to be a vertex of the $H$-graph).

\item[{$\cE_4$:} ] An $\cE_3$ edge such that for any $\hat{S}$-vertex $x \in
  e$, $1-\theta_4$ fraction of the edges reverse branching from $x$
  are in $\cE_1$.

\end{itemize}
\end{definition}

\noindent{\bf Remark:} There $\cE_2$ edges, they will be defined
below in \S\ref{subsection:circuits}.

Choose $1 \GG  \nu_3 > \nu_2 > 0$. The $\nu$'s will be functions
of $\delta_0$ which tend to $0$ as $\delta_0 \to 0$.

\begin{definition}[Good Vertices]
\label{def:good:vertices}
The following defines sets of ``good'' vertices. See also
Definition~\ref{def:good:edges}.
\begin{itemize}
\item[{$\cV_1$:} ] The set of ``good'' vertices as defined in the
  previous section.
\item[{$\cV_2$:} ] In $\cV_1$ and $1 -\nu_2$ fraction of the
outgoing edges are in $\cE_1$.

\item[{$\cV_3$:} ] In $\cV_2$ and $1-\nu_3$ fraction of the
outgoing
  edges are in $\cE_4$.

\item[{$\cV_4$:} ] In $\cV_3$ and is not a strange vertex (see
  Definition~\ref{def:strange:vertex} below).

\end{itemize}
\end{definition}

\begin{lemma}
\label{lemma:most:V1}
If $L' \GG L$, we can choose a horocycle $H$ such that for the $H$-graph
$\cG(H)$, $1-\delta_1$ fraction of vertices are in $\cV_1$. Here,
$\delta_1$ is a function of $\delta_0$ which tends to $0$ as $\delta_0
\to 0$.
\end{lemma}

\bold{Proof.}
Note that $\phi^{-1}(B(L'))$ has small boundary area (compared to the
volume). Now tile $\phi^{-1}(B(L'))$ by boxes $B(L)$. Since $L' \GG
L$, most boxes are completely in the interior of $\phi^{-1}(B(L'))$.

Let $\cU$ denote the set where we know the map is locally standard (but could
be right side up or upside down).
Note that for every box $B(L)$, $|\cU \cap B(L)| \ge 0.999 |B(L)|$.

This implies that for most $H$,
\begin{displaymath}
| H \cap \phi^{-1}(B(L')) \cap \cU | \ge 0.99 | H \cap
\phi^{-1}(B(L')) |
\end{displaymath}
Then for such $H$, $\cV_1$, which consists of vertices on $\phi(H)
\cap B(L') \cap \phi(\cU)$, satisfies the conditions of the lemma.
\qed\medskip

We now fix $H$ such that Lemma~\ref{lemma:most:V1} holds.

\begin{lemma}
\label{lemma:most:E1} At least $1-\epsilon_1$ fraction of the
edges of $\cG(H)$ are in $\cE_1$. Here, $\epsilon_1$ is a function
of $\delta_0$ which tends to $0$ as $\delta_0 \to 0$.
\end{lemma}

\bold{Proof.} Recall that $m\leq M(v) \leq M$ where $M(v)$ is the
degree of $v$. This implies that

\begin{displaymath}
\label{equation:ratio:edges:vertices} \frac{1}{m} \leq
\frac{|\cV(H)|}{|\cE(H)|} \leq \frac{1}{M}.
\end{displaymath}

Since each edge not in $\cE_1$ must be quasi-incident on a vertex
not in $\cV_1$ and each vertex is incident to at most $M$ edges,
we have:

\begin{displaymath}
| \cE_1^c | \le 2 M | \cV_1^c|.
\end{displaymath}

Combined with equation $(\ref{equation:ratio:edges:vertices})$
this implies
\begin{displaymath}
\frac{|\cE_1^c|}{|\cE(H)|} \le 2\frac{M}{m}
\frac{|\cV_1^c|}{|\cV(H)|}
\end{displaymath}
Thus the lemma follows from Lemma~\ref{lemma:most:V1}.
\qed\medskip

\begin{lemma}
\label{lemma:most:V2}
At least $1-\delta_2$ fraction of the vertices of $\cG(H)$ are in $\cV_2$.
Here, $\delta_2$ is a function of $\delta_0$ which tends to $0$ as $\delta_0
\to 0$.
\end{lemma}

\bold{Proof.} This follows immediately from
Lemma~\ref{lemma:most:E1}. \qed\medskip

\begin{lemma}
\label{lemma:most:E3}
At least $1-\epsilon_3$ fraction of the edges of $\cG(H)$ are in $\cE_3$.
Here, $\epsilon_3$ is a function of $\delta_0$ which tends to $0$ as $\delta_0
\to 0$.
\end{lemma}

\bold{Proof.} In view of Lemma~\ref{lemma:most:V2}, it enough to
prove that for any $v \in \cV_2$, almost all the edges outgoing from
$v$ belong to $\cE_3$.

Suppose $v \in \cV_2$. Let $\cE(v)$ denote all the edges which are
incident to $v$. We know that most edges in $\cE(v)$ belong to
$\cE_2$, i.e.
\begin{equation}
\label{eq:E1:Ev:small} | \cE_2^c \cap \cE(v) | \le \delta_2 |
\cE(v)|
\end{equation}
Let $A_v = \cE(v) \cap \cE_3^c$ denote the edges outgoing from $v$
which are not in $\cE_3$. We know that for any $e \in A_v$, there
exists $x \in e$ such that at least $\theta_3$ of the edges
branching from $e$ at $x$ are not in $\cE_1$. Thus there exists a
neighborhood $U$ of $e$ such that
\begin{displaymath}
|\cE_2^c \cap U \cap \cE(v)| \ge \theta_3 | U \cap \cE(v)|
\end{displaymath}
We thus get a cover of $A_v$ by $U$'s. Then by Vitali's covering
lemma, there exists disjoint $U_j$ such that
\begin{displaymath}
\sum_{j=1} |U_j| \ge \frac{1}{2} |A_v|.
\end{displaymath}
Thus,
\begin{displaymath}
|A_v| \le 2 \sum |U_j| \le \frac{2}{\theta_3} \sum | U_j \cap
\cE_2^c \cap \cE(v)| \le \frac{2}{\theta_3} | \cE_2^c \cap \cE(v) |
\end{displaymath}
Then, by (\ref{eq:E1:Ev:small}),
\begin{displaymath}
|A_v| \le \frac{2 \delta_2}{\theta_3} | \cE(v)|.
\end{displaymath}
We now choose $\theta_3 = \epsilon_3 = \sqrt{2 \delta_2}$.
\qed\medskip

\begin{lemma}
\label{lemma:most:E4} At least $1-\epsilon_4$ fraction of the edges
of $\cG(H)$ are in $\cE_4$. Here, $\epsilon_4$ is a function of
$\delta_0$ which tends to $0$ as $\delta_0 \to 0$.
\end{lemma}

\bold{Proof.}
It follows immediately from Lemma~\ref{lemma:most:E3} that $1 - 2
\epsilon_3$ proportion of the for non-leaf edges have the reverse
branching property.

Let $\alpha = \epsilon_1^{1/6}$. If the proportion of the leaf
edges is at most $\alpha$, we are already done (with $\epsilon_4 =
2 \epsilon_3 + \alpha$). Thus we may assume that the proportion of
leaf edges is at least $\alpha$.

Let $Y$ be the set of all vertical paths in $B(L)$ going from top to
bottom, and $Y' \subset Y$ is the subset consisting of paths which
pass through a vertex not in $\cV_1$. Let $D(\gamma)=1$ if
$\gamma{\in}Y'$ and $D(\gamma)=0$ otherwise.
By lemma \ref{lemma:most:E1}, we have:
\begin{equation}
\label{eq:temp} \sum_{\gamma \in Y} D(\gamma) \le {2 \epsilon_1}
|\cE(H)|.
\end{equation}
\relax From (\ref{eq:temp}),
\begin{displaymath}
\sum_{\gamma \in Y} D(\gamma) \le \epsilon_1 |\cE(H)| \le \frac{
\epsilon_1}{\alpha} |\cE_{leaf}|,
\end{displaymath}
where $\cE_{leaf} \subset \cE(H)$ denotes the set of leaf edges.
For a point $v \in \partial B(R)$, let $Y_v$ denote the set of
geodesics emanating from $v$. We get
\begin{displaymath}
\sum_{v \in \partial B(R)} \sum_{\gamma \in Y_v} D(\gamma) \le
\sum_{v \in \partial B(R)} \frac{\epsilon_1}{\alpha}
|\cE_{leaf}(v)|,
\end{displaymath}
where $\cE_{leaf}(v)$ denotes the set of leaf edges emanating from
$v$. Let $\theta' = \epsilon_1^{2/3}$, and let
\begin{displaymath}
P = \left\{ v \in \partial B(R) \st \sum_{\gamma \in Y_v} D(\gamma) >
\theta' |\cE_{leaf}(v)| \right\}
\end{displaymath}
Note that
\begin{displaymath}
\sum_{v \in P} |\cE_{leaf}(v)| \le \sum_{v \in P} \frac{1}{\theta'}
\sum_{\gamma \in Y_v} D(\gamma) \le \frac{1}{\theta'} \sum_{v \in Y}
D(\gamma) \le \frac{\epsilon_1}{\alpha \theta'} |\cE_{leaf}|
\end{displaymath}
Thus, since we choose $\alpha$ and $\theta'$ so that
$\frac{\epsilon_1}{\alpha \theta'} \ll 1$, it is enough to prove
that for $v \not\in P$, most of the edges in $\cE_{leaf}(v)$ are in
$\cE_4$.

Now assume $v \not\in P$. Thus we have
\begin{displaymath}
\sum_{\gamma \in Y_v} D(\gamma) <
\theta' |\cE_{leaf}(v)|
\end{displaymath}

Choose $\theta_4 = \epsilon_1^{1/12}$. Let $A_v = \cE_{leaf}(v) \cap
\cE_4^c$ denote the leaf edges outgoing from $v$ which are not in
$\cE_4$. We know that for any $e \in A_v$, there exists $x \in e$
such that at least $\theta_4$ of the edges branching from $e$ at $x$
are not in $\cE_2$. Thus there exists a neighborhood $U \subset Y_v$
with $e \in U$ such that
\begin{displaymath}
|\cE_2^c \cap U| \ge \theta_4 | U |
\end{displaymath}
hence using the definition of $\cE_2$,
\begin{displaymath}
\sum_{\gamma \in U} D(\gamma) \ge \theta_2 \theta_4 |U|.
\end{displaymath}
We thus get a cover of $A_v$ by $U$'s. Then by Vitali, there
exists disjoint $U_j$ such that
\begin{displaymath}
\sum_{j=1} |U_j| \ge \frac{1}{2} |A_v|.
\end{displaymath}
Thus,
\begin{displaymath}
|A_v| \le 2 \sum_j |U_j| \le \frac{2}{\theta_2 \theta_4} \sum_j
\sum_{\gamma \in   U_j} D(\gamma) \le \frac{2}{\theta_2 \theta_4}
\sum_{ \gamma \in Y_v} D(\gamma) \le \frac{2 \theta'}{\theta_2 \theta_4}
|\cE_{leaf}(v)|.
\end{displaymath}
Since $\theta_2 = \epsilon_1^{1/2}$, $\frac{2 \theta'}{\theta_2
  \theta_4}=\epsilon_1^{\frac{1}{12}}$ and the lemma follows.
\qed\medskip

\begin{lemma}
\label{lemma:most:V4} At least $1-\delta_4$ fraction of the
vertices of $\cG(H)$ are in $\cV_3$. Here, $\delta_4$ is a
function of $\delta_0$ which tends to $0$ as $\delta_0 \to 0$.
\end{lemma}

\bold{Proof.} This follows from Lemma~\ref{lemma:most:E4}.
\qed\medskip

Let $H_*$ be an horocycle intersecting $B(L')$. We say that an
$S$-vertex on $w$ on $H_*$ is {\em marked} by a $\cV_1$ $H$-vertex
$v$ if the cloud of $v$ contains a point of $H_*$, and also $h(v) =
h(H_*) + O(\rho_2)$, and also the coordinates of $v$ and $w$ along
$H_*$ must agree up to $O(\rho_2)$. (In particular the orientation
of $v$ must be such that the coordinate of $v$ along $H_*$ is not
``fuzzy'').

% Given an $H$-vertex $v$ we can associate to it a vertex in the
% range.  For bad vertices, this is obvious.  For a good vertex, we
% take $\hat \phi(v)$ where $v$ is the vertex in the domain where
% edges in the equivalence class of $v$ come close. Given a vertex
% with representative $v$ and a horocycle $H$ through $v$, we say
% that another vertex $w$ is {\em along $H$} if either

% \begin{enumerate}
% \item $w$ is a bad vertex and $w$ is on $H$ or \item $w$ is a good
% vertex and $w$ is within $\epsilon R$ of $H$.
% \end{enumerate}

\begin{definition}[Strange Vertex]
\label{def:strange:vertex} An $H$-vertex $v \in \cV_3$ is called
{\em strange} if there is an  horizontal segment (i.e. piece of
horocycle) $K$ marked by $v$ such that more then $1-\nu_4$ fraction
of the $S$-vertices on $K$ are marked by $H$-vertices which are
$\cV_1$ but not in $\cV_3$.
\end{definition}

% Let $A(v,H)$ be the set of vertices along $H$.  Then there is a
% projection $\pi: A(v,H) \rightarrow H$. For bad vertices, this
% projection is the identity.  For a good vertex $w$, we consider
% the actual vertex in the domain representing $w$ project it to
% $\tilde w$ on the domain horocycle mapped by $\hat \phi$ to $H$
% and take $\hat \phi(\tilde w)$. This map is one to one since if
% two vertices have the same image then they are easily seen to
% belong to the same equivalence class. \mc{Say this better.}

\begin{lemma}
\label{lemma:most:V6} At least $1-\delta_6$ fraction of the
vertices of $\cG(H)$ are in $\cV_4$ (i.e are in $\cV_3$ and not
strange). Here, $\delta_6$ is a function of $\delta_0$ which tends
to $0$ as $\delta_0 \to 0$.
\end{lemma}

\bold{Proof.} Let $v_1, \dots, v_m$ be the strange vertices, and let
$K_1, \dots, K_m$ be horocycle segments marked by the strange
vertices. The $K_i$ are not quite uniquely defined, but we address
this issue below.

Note that the number of $H$-vertices which can mark a given
$S$-vertex is $O(\rho_2)$. Indeed, any two such vertices must be
within $O(\epsilon' R)$ of each other, which means that they must
have come from the same good box, which implies that heights and
their transverse coordinates must agree. (Recall that the vertices
which come from near the edges of a good box are automatically not
in $\cV_1$).

The same argument shows that one can choose the horocycle segments
$K_i$ so that for $i
\ne j$, $d(K_i, K_j) > 3 D \epsilon' R$. Now we can apply the Vitali
covering lemma to the $K_i$. This lemma applies since each $K_i$ is
one-dimensional and the different $K_i$ do not interact with each
other. Also the density of the $\cV_1$ vertices which are not in
$\cV_3$ is small by Lemma~\ref{lemma:most:V4}. This implies that the
strange vertices are a small fraction of all the vertices.
\qed\medskip

\subsection{Circuits}
\label{subsection:circuits}

\bold{The projection $\pi_H$ and the function $\rho_H( \cdot,
\cdot)$.} Let $H$ be a horocycle. Let $\pi_H: \Sol \to \half^2$
denote the orthogonal projection to the hyperplane orthogonal to
$H$. We let $\rho_H(p,q)=(\pi_H(p)|\pi_H(q))_{\pi_H(H)}$ be the
Gromov product of $\pi_H(p)$ and $\pi_H(q)$ with respect to
$\pi_H(H)$ in $\half^2$. Recall that for three points $x,y,z$ in a
metric space $X$, the Gromov product is defined as:

$$(y|z)_x=\frac{1}{2}\{d_X(x,y)+d_X(x,z)-d_X(z,y)\}.$$

\noindent Let $\gamma_{yz}$ be the geodesic joining $y$ to $z$. In a
$\delta$-hyperbolic space $X$ satisfies

$$d_X(\gamma_{yz},x)-\delta \leq (y|z)_x \leq d_X(\gamma_{yz},x)$$

\noindent see e.g. \cite[Lemma 2.17]{GhdlH}. We note the following
properties of $\rho_H$:

\begin{lemma}
\label{lemma:properties:rho:H}
\begin{itemize}

\item[{\rm (i)}] Suppose $d(p',p) \ll d(p,H)$, $d(q',q) \ll
d(q,H)$,
  and $\rho_H(p,q) \ll \min(d(p,H),d(q,H))$. Then,
\begin{displaymath}
\rho_H(p,q) \approx \rho_H(p',q').
\end{displaymath}

\item[{\rm (ii)}] Suppose $h(p') < h(p)$, $h(q') < h(q)$, the
points
  $p$ and $p'$ can be connected by a vertical geodesic, and the same
  for the points $q$ and $q'$. Suppose also $d(p,H) \GG \rho_H(p,q)$
  and $d(q,H) \GG \rho_H(p,q)$. Then,
\begin{displaymath}
\rho_H(p,q) \approx \rho_H(p',q').
\end{displaymath}

\item[{\rm (iii)}] If $\rho_H(p,q) > s$ and $\rho_H(q,q') > s$,
then $\rho_H(p,q') > s$ (up to a small error).
\end{itemize}
\end{lemma}

\bold{Proof.} The statements (i), (ii) and (iii) are standard
hyperbolic geometry. In particular (iii) follows immediately from
the ``thin triangle'' property.
\qed\medskip

In the following lemma, the horocycle $H$ is assumed to be an $y$
horocycle.  An analogous lemma, with a few sign changes, holds for
$x$ horocycles.

\begin{lemma}
\label{lemma:nodrift:one} Suppose $p, q \in X(n)$ are connected by
a path $\hat{\gamma}$ in such that
\begin{equation}
\label{eq:nodfift:low} h(x) \le h(H) -\rho_4 \text{ for all $x \in
\gamma$,}
\end{equation}
Further assume the initial segments of $\gamma$ at both $p$ and
$q$ are vertical geodesics going down for length at least
$\epsilon' R$, that $\gamma$ stays below $h(H)-R'$ except on these
initial segments and that the length of $\gamma$ is less than
$e^{\epsilon'
  R}$.  Then, $\rho_H(p,q) > \Omega(\rho_4)$.
\end{lemma}

\bold{Proof.} This is standard hyperbolic geometry applied to
$\pi_H(\gamma)$.
\qed\medskip

\bold{Notation.} An $\cE_2$ edge is a monotone vertical path in the
$\hat{S}$-graph which is a subset of an $\cE_1$ edge (or possibly
a subset the extension of an $\cE_1$ edge by at most $\rho_4$ at
each end).

% We introduce a new class of paths in $\Sol$ which are not quite
% subsets of the $H$-graph.  An $\cE_2$ edge is a monotone geodesic
% segment in $\Sol$ obtained by either cutting or extending an
% $\cE_1$ edge by length at most $\rho_4$ on each end.  By the
% definition of the $H$ graph and good vertices, it follows that the
% ends of an $\cE_2$ edge are in good boxes and that they are
% entirely contained in the subset where $\phi$ is close to $\hat
% \phi$.

\begin{lemma}
\label{lemma:nodrift:two} Suppose $\gamma = \overline{p_0 q_0}$ is
an $\cE_2$ edge going up from an $x$-oriented irregular
$\hat{S}$-vertex (or going down from an $y$-oriented irregular
$\hat{S}$-vertex). Suppose $p \in \gamma$ is within the same
$B_i(R)$ as $p_0$, and $q \in \gamma$ is within the same
$B_{i'}(R)$ as $q_0$ and $d(p,p_0)$ and $d(q,q_0)$ is at least $10
\epsilon' R$. Then the following hold:
\begin{itemize}
\item[{\rm (i)}] Except at near it endpoints, $\gamma$ never
passes through any irregular $\hat{S}$-vertices. In other words,
in its' interior, $\gamma$ never comes within $R'$ of a good $H$
vertex. \item[{\rm (ii)}] We have
  $\rho_H(\hat{\phi}^{-1}(p),\hat{\phi}^{-1}(q)) > \Omega(\rho_4)$.
\end{itemize}
\end{lemma}

\bold{Remark.}
In the above $\rho_H(\hat{\phi}^{-1}(p), \hat{\phi}^{-1}(q))$ is
well defined since for an $\hat{S}$-vertex $v$,
$\pi_H(\hat{\phi}^{-1}(v))$ is well defined (even though
$\hat{\phi}^{-1}(v)$ may not be).

\bold{Informal outline of proof.} We first outline the proof, and
then give the full argument. Consider $\phi^{-1}(\gamma)$. Note
that below height $h(H)-R'$, $\phi^{-1}(\gamma)$ cannot move
transverse to $H$ because it is of length at most $O(L)$.  Because
of this, whenever $\gamma$ attempts to cross above height
$h(H)-R'$ it must does so in the image of $W(H)$. Consider the
point $q'$ where it does so. Since $\gamma$ cannot hit a bad
shadow vertex, $q'$ must be essentially in $U \cap B_i(R) \cap
W(H)$. But then, by the definition of the $\hat{S}$-graph,
$\gamma$ must hit an $H$-vertex. Thus, $q'$ is near the endpoint
of $\gamma$, and thus (i) holds. Now (ii) follows from
Lemma~\ref{lemma:nodrift:one} since we know that
$\phi^{-1}(\gamma)$ has not passed above height $h(H)-R'$ except
near the endpoints.

\bold{Proof.} Let $p_1$ be the first place where $\gamma$ hits
$\hat{W}(H)$. Then, since $\gamma$ cannot hit a bad shadow vertex,
there exists $p_1' \in U_* \cap W(H)$ such that $\hat
\phi(p_1')=p_1$ and $d(\phi^{-1}(p_1),\hat{\phi}^{-1}(p_1)) =
O(\epsilon' R)$. Note that $p_1'$ and $\phi{\inv}(p_1)$ are both
$\Omega(\epsilon' R)$ from the sides of $W(H)$.

Let $p'_2$ be the next point after $p_1'$ when $\phi^{-1}(\gamma)$
intersects $\tilde{B} \cap \{x \st h(x) = h(H) - R' \}$ at
$\phi{\inv}(p_2)$.  Since $\gamma$ is an $\cE_2$ edge and in
particular a vertical geodesic, we know $d(p'_1,p'_2)$ is
$\Omega(\beta R)$. By the choice of $p'_2$ and the definition of
$\cE_2$ edge, $\phi{\inv}\gamma$ never hits a shadow vertex
between $p'_1$ and $p'_2$. This fact and the fact that
$|\phi^{-1}(\gamma)|<O(L)$ imply that $p'_2$ must be in $W(H)$.
Since $\gamma$ is $\cE_2$, $p'_2$ is not a bad shadow vertex and
in particular is away from the edge of $W(H)$. Together this
implies that $p'_2$ is in $W(H){\cap}\tilde B$ and that the
continuation of $\gamma$ past $p_2=\hat \phi(p'_2)$ must, by the
definition of the $\hat{S}$ and $H$ graphs, hit an $H$-vertex.
Since $\gamma$ is an $\cE_2$ edge and does not contain good
vertices in it's interior, this implies that $p_2$ and $q_0$ are
in the same box and that the segment from $p_2$ to $q_0$ contains
$q$. Now by Lemma~\ref{lemma:properties:rho:H},
\begin{displaymath}
\rho_H(\hat{\phi}^{-1}(p),\hat{\phi}^{-1}(q)) =
\rho_H(\hat{\phi}^{-1}(p_1), \hat{\phi}^{-1}(p_2)) \approx
\rho_H(\phi^{-1}(p_1), \phi^{-1}(p_2)) \ge \Omega(\rho_4)
\end{displaymath}
\qed\medskip

\begin{lemma}
\label{lemma:comb:sameH} Suppose $\overline{p_0 q_0}$ is an
$\cE_2$ edge (which goes up from a $x$-oriented vertex and down
from an $y$-oriented vertex),  $\rho_3 \GG s \GG \rho_1$, and $p$
(resp. $q$) is on $\gamma$ distance $s$ away from $p_0$ (resp.
from $q_0$). Then there exists a horocycle $H'$ such that $p$ and
$q$ are within $O(\rho_1)$ of $\hat{\phi}(H')$.
\end{lemma}

\bold{Proof.} Choose points $p'$ and $q'$ on $\overline{p_0q_0}$
close to where $\overline{p_0q_0}$ enters the respective good
boxes. Applying Lemma~\ref{lemma:nodrift:two} we see that
$\rho_H(\hat{\phi}^{-1}(p'), \hat{\phi}^{-1}(q')) >
\Omega(\rho_4)$.
%However, $d(\hat{\phi}^{-1}(p), H) \ll \rho_3 < \rho_4$.
By the usual $\delta$-thin triangle properties, this implies that
the geodesic segments
$\overline{\pi_H(\hat{\phi}^{-1}(p))\pi_H(H)}$ and
$\overline{\pi_H(\hat{\phi}^{-1}(q))\pi_H(H)}$ stay close till
roughly for roughly $\rho_4$  units from $\pi_H(H)$. Since
$d(\hat{\phi}^{-1}(p), H) \ll \rho_3 < \rho_4$ and similarly for
$\hat{\phi}^{-1}(q)$ this implies that $\pi_H(\hat{\phi}^{-1}(p))$
and $\pi_H(\hat{\phi}^{-1}(q))$ are within $2\delta$ of the same
vertical geodesic through the point $\pi_H(H)$. But since they are
at the same height, this implies that $\pi_H(\hat{\phi}^{-1}(p)) =
\pi_H(\hat{\phi}^{-1}(q))$. \qed\medskip

Suppose $H'$ is a horocycle obtained by moving
up less than $\rho_3$ from $H$.
Recall that the set $\hat{\phi}(H')$ is a well defined
subset of the $\hat{S}$-graph (see \S\ref{subsection:build:H:graph}).
We always assume that $\hat{\phi}(H')$ runs along
vertices in the $\hat{S}$-graph (or else project it). By Lemma
\ref{lemma:comb:sameH}, given any collection of $\cE_2$ edges with
(some) endpoints on $H$, we may replace them with $\cE_2$ edges
with (some) endpoints on $H'$.

\begin{lemma}[Illegal Circuit]
\label{lemma:illegal:1:circuit} Suppose $n$ is some finite even
integer which is not too large (we will use $n=4$ and $n=6$), and
for $0 \le i \le n-1$, $p_i$ are $\hat{S}$-vertices. Also suppose that
for $0 \le i \le n-1$, $ \overline{p_{i-1} p_i}$ are subsets of
$\cE_2$ edges, where $i-1$ is considered mod $n$. For $1 \le i \le
n$, let $r^\pm(p_i)$ denote the maximum distance the geodesic
$\overline{p_{i \pm 1} p_i}$ can be continued beyond $p_i$ while
remaining a subset of an $\cE_2$ edge, and let $r(p_i) = \max(
r^+(p_i),r^-(p_i))$.

Suppose there is an index $k$ such that $r(p_k) \ll \rho_4$, and
for all $i \ne k$, $r(p_i) > r(p_k)+2{\rho_1}$. Then
$\overline{p_{k-1} p_k}$ and $\overline{p_k p_{k+1}}$ cannot have
only the point $p_k$ in common.
\end{lemma}

\noindent {\bf Remark:} Roughly, the point of the lemma is that
one cannot find a loop of length $O(L)$ through a point on the
horocycle which begins by going up in two distinct directions
unless the loop comes back to the original horocycle.

\bold{Proof.} Without loss of generality, $k=0$. Let $H'$ be the
horocycle passing thorough $\hat{\phi}^{-1}(p_0)$. By Lemma
\ref{lemma:comb:sameH} and the discussion following, we can
consider $H'$ in place of $H$, namely we can replace all $H$
vertices that occur in our arguments with vertices in $H'$. Let
$p_{i-1}^+$ be the first time when $\overline{p_{i-1} p_i}$ leaves
$\hat{\phi}(\tilde{B} \cap W(H))$, and let $p_{i-1}^-$ be the last
time when $\overline{p_{i-1} p_i}$ enters $\hat{\phi}(\tilde{B}
\cap W(H))$ (so $d(p_{i-1},p_{i-1}^+) \approx R' \le R$, and
$d(p_i^-, p_i) \approx R' \le R$). By applying
Lemma~\ref{lemma:nodrift:two} to each segment
$\overline{p_{i-1}^+,
  p_i^-}$ we see that $\rho_{H'}(\hat{\phi}^{-1}(p_{i-1}^+),
\hat{\phi}^{-1}(p_i^-)) \ge \Omega(\rho_4)$.

Now, by assumption, for all $i \in [0,n-1]$ except $k=0$,
$\rho_{H'}(\hat{\phi}^{-1}(p_i^-),\hat{\phi}^{-1}(p_i^+)) \ge
2{\rho_1}$, but for $i=0$,
$\rho_{H'}(\hat{\phi}^{-1}(p_0^-),\hat{\phi}^{-1}(p_0^+)) \le
\rho_1$. This is a contradiction to
Lemma~\ref{lemma:properties:rho:H} (iii). \qed\medskip

\subsection{Families of geodesics}
\label{section:families}

Let $B[\lambda]$ be a box in $X(n')$
of combinatorial size $\lambda$ (i.e. the number of
edges from the top to the bottom is $\lambda$, and the distance from
the top to the bottom is $\rho_1 \lambda$). Let $\gb$ be the branching
number of each vertex (i.e the valence of each vertex
counting both up and down branching is $2 \gb$).
Note that $\gb$ is related to the the constant $B_X'$ of
\S\ref{subsection:shadows} by $\gb^{2\lambda} = e^{B_X' \rho_1
  \lambda}$, so $\log \gb = B_X' \rho_1/2$.

Thus the number of $\hat{S}$-vertices on the top edge of
$B[\lambda]$ is $\gb^\lambda$, and so is the number of
$\hat{S}$-vertices on the bottom edge. The total number of vertical
geodesics in $B[\lambda]$ is $\gb^{2 \lambda}$.

\begin{lemma}
\label{lemma:number:vertices} The number of vertices of the $H$-graph
in $B[\lambda]$ is at most $c_9(\rho_1) \gb^\lambda$.
\end{lemma}

\bold{Proof.} Apply Lemma~\ref{lemma:trapping} with $Q$ the union of
the top edge and the bottom edge of the box.
\qed\medskip

Given a box $B(D)$ and a vertical geodesic segment $\gamma$ of
length $D$ in $B$, we say $\gamma$ {\em is through} if $\gamma$
does not hit any $H$ vertex in $B$.  The following lemma applies
to families of geodesics in a box.  Note that the geodesics are
not assumed to be part of the $H$-graph. The point of the lemma is
that if too many paths through the box are blocked by good
vertices, then some good vertex must block many paths.  This
really only depends on the fact that there are not too many good
vertices in the box.

\begin{lemma}
\label{lemma:family} Let $B[\lambda]$ be a box of combinatorial size
$\lambda$. Suppose $\cF$ is a family of vertical geodesics (actually
monotone paths in the modified $\hat{S}$-graph going from the top of
$B[\lambda]$ to the bottom) with the following properties:
\begin{itemize}
\item[{\rm (a)}] Each geodesic in $\cF$ does not hit any bad vertices
 % (i.e. is a union of edges in $\cE_1$.)

\item[{\rm (b)}] $|\cF|$ (i.e. the number of geodesics in $\cF$) is at
  least $\sigma \gb^{2 \lambda}$, where $0 < \sigma < 1$.

%\item[{\rm (c)}] There exists $\rho \in \natls$ (we will always
%use
%  $\rho = \rho_2$) such that for any
%  vertex $v \in \cV_1$ not on the bottom edge or within $\rho$ of the
%  top edge of $B[\lambda]$,
%  any two geodesics in $\cF$ which pass through $v$ stay together
%  for at least $\rho$ $S$-edges.
%\end{itemize}

%Then, all but $\frac{c_9(\rho_1)}{\beta^{\rho}} \beta^{2\lambda}$
%of the geodesics in $\cF$ are through (i.e. do not contain any
%$H$-vertices in their interior). Thus, if $\sigma \GG
%\frac{c_9(\rho_1)}{\beta^{\rho}}$, almost all of the geodesics in
%$\cF$ are through.
%\end{lemma}

\item[{\rm (c)}] For some $\rho \in \natls$ (we will always use
  $\rho = \rho_2$), fewer than $1-\frac{c_9(\rho_1)}{\gb^{\rho}}
\gb^{2\lambda}$ of the geodesics in $\cF$ are through (i.e. do not
contain any $H$-vertices in their interior).
\end{itemize}

Then there exists a
  vertex $v \in \cV_1$ not on the bottom edge or within $\rho$ of the
  top edge of $B[\lambda]$, and
  two geodesics in $\cF$ which pass through $v$ and stay together
  for fewer than $\rho$ $\hat{S}$-edges.
\end{lemma}

Thus, if $\sigma \GG \frac{c_9(\rho_1)}{\gb^{\rho}}$, almost all
of the geodesics in $\cF$ are through unless we have a
configuration as described in the conclusion of the lemma.

\bold{Proof.} Let $\cF_0$ denote the family of all vertical geodesics
on the unmodified $\hat{S}$-graph,
passing from the top of $B[\lambda]$ to the bottom. Clearly $|\cF_0| =
\gb^{2 \lambda}$. Note that $B[\lambda]$ has $\lambda \gb^{\lambda}$
$\hat{S}$-vertices, and each geodesic contains $\lambda$ $\hat{S}$-vertices.
This implies that each $\hat{S}$-vertex lies on $M = \gb^\lambda$
geodesics in $\cF_0$.

Now suppose $v$ is an $H$-vertex in $\cV_1$, and $v$ is not on the
bottom edge or within $\rho$ of the top edge. Assuming the
conclusion of the lemma fails then $v$ can belong to at most $M
\gb^{-\rho}$ geodesics in $\cF$. Thus, using
Lemma~\ref{lemma:number:vertices}, we see that the total number of
geodesics in $\cF$ which pass through a vertex in $\cV_1$ in $B$
is at most
\begin{displaymath}
M \gb^{-\rho} c_9(\rho_1) \gb^{\lambda} = \frac{c_9(\rho_1)}{
  \gb^{\rho}} \gb^{2 \lambda}
\end{displaymath}
which implies that all but $\frac{c_9(\rho_1)}{
  \gb^{\rho}} \gb^{2 \lambda}$ of the geodesics are through,
  contradicting $(c)$.

\qed\medskip

\bold{Convention.} For the remainder of this subsection, we assume
that we have a $x$ horocycle $H$ whose image (at least in some
initial box) is $x$-oriented.  The proof proceeds by extending the
set on which the image is horocycle and so all points in the $H$
graph we consider will be $x$-oriented.

\bold{The intervals $I_\lambda(v)$ and $I_\lambda'(v)$.} For any
$\hat{S}$-vertex $v$, let $I_{\lambda}(v)$ denote the set of
vertices on the same $x$-horocycle as $v$ which are within
combinatorial distance $2 \lambda$. Let $I_\lambda'(v)$ denote the
set of vertices which can be reached from $v$ by a monotone path
going up for exactly $\lambda$ steps (so $I_\lambda'(v)$ is a
piece of $y$-horocycle). Note that for $DL$ graphs, each point of
$I_\lambda(v)$ is connected to each point of $I_\lambda'(v)$ by a
monotone path of length $\lambda$, and for any $w \in
I_\lambda(v)$, $I_\lambda(w) = I_\lambda(v)$ and $I_\lambda'(w) =
I_\lambda'(v)$.  For $\Sol$ slightly more complicated variants of
these statements hold, for instance for any $w \in I_\lambda(v)$,
$I_\lambda(w)$ and $I_\lambda(v)$ intersect in a set that contains
more than half the measure of each.  And the the relative measure
of this intersection in each set is close to one, unless $w$ is
close to an edge of $I_\lambda(v)$.

Let $v$ be any $\hat{S}$-vertex. Let ${\mathcal U}(v,\lambda)$ denote
the set of distinct monotone geodesic segments in the $\hat{S}$ graph
going up from $v$ for distance exactly $\lambda$. Then
$U(v,\lambda)$ is the is the set of geodesics joining $v$ to
$I_{\lambda}'(v)$. Similarly, we let ${\mathcal D}(w,\lambda)$ be
the set of distinct monotone geodesics segments in the $\hat{S}$-graph
going down distance $\lambda$ from $w$.  If
$w{\in}I_{\lambda}'(v)$ then $\mathcal D(w,\lambda)$ is the set of
monotone geodesics joining $w$ to points in $I_{\lambda}(v)$

\begin{proposition}[Extension of Horocycles I]
\label{prop:extension} Suppose $v \in \cV_3$. Suppose $\sigma \GG
\eta \GG c_9(\rho_1)/\gb^{\rho_2}$ and suppose $\lambda$ is such
that at least $\sigma$-fraction of the edges going up from $v$ are
$\cE_4$ edges of length at least $\lambda+\rho_2$. Then at least
$1-O(\eta)$ fraction of the $\hat{S}$-vertices in $I_\lambda(v)$
are in fact $H$-vertices.
\end{proposition}

\bold{Proof.} We assume that $v \in \cV_3$ and $\hat{\phi}(H)$ is
oriented as a $x$-horocycle near $v$. Let $E$ denote the set of
$\cE_4$ edges coming out of $v$ which have length at least
$\lambda+\rho_2$. Let $E_{\lambda}$ be the set of vertices in
$I'_\lambda(v)$ which are on of $\lambda+\rho_2$ unobstructed
geodesics leaving $v$. By assumption, we have
\begin{equation}
\label{eq:E:lambda:j:large} |E_{\lambda}| \ge \sigma
\gb^{\lambda}.
\end{equation}
We now let $\cF_0'=\bigcup_{w{\in}E_{\lambda}}\cD_{\lambda}(v)$
and let $\cF'$ be all the geodesics segments in $\cF'_0$ which do
not contain a bad vertex. Assume for a contradiction that many
geodesics in $\cF'$ are not through, i.e. that  $(c)$ of
Lemma~\ref{lemma:family} holds for $\cF'$. We verify that
Lemma~\ref{lemma:family}$(a)$ and $(b)$ hold for $\cF'$. Since $v
\in \cV_3$,
\begin{equation}
\label{eq:bounds:Fj:prime} (1-\theta_4) |E_{\lambda}|
\gb^{\lambda}   \le |\cF'|
   \le |E_{\lambda}| \gb^{\lambda}
\end{equation}
Note that by (\ref{eq:E:lambda:j:large}) and
(\ref{eq:bounds:Fj:prime}), we have $|\cF'| \ge \sigma
\gb^{2(\lambda)}$. Hence Lemma~\ref{lemma:family} (b) holds.  Note
that all the geodesics in $\cF'$ end at points of $I_\lambda(v)$.

Now by Lemma \ref{lemma:family} there exists $w \in \cV_1$ with
$h(w)
> h(v)$ and an $\hat{S}$-vertex $w_1$ with $h(w_1) > h(w)$ and $d(w,w_1)
< \rho_2$ so that at least two geodesics in $\cF'$ meet at $w_1$
and continue to $w$. (See figure~\ref{fig:extension:1}). Let $x
\in I_\lambda'(v)$ and $y \in I_\lambda'(v)$ be the starting
points of these two geodesics.

\makefig
{Proof of Proposition~\ref{prop:extension}. \\
The filled boxes denote $H$-vertices.}{fig:extension:1}{
\begin{picture}(0,0)%
\includegraphics{circuit1.pstex}%
\end{picture}%
\setlength{\unitlength}{3947sp}%
\begingroup\makeatletter\ifx\SetFigFont\undefined%
\gdef\SetFigFont#1#2#3#4#5{%
  \reset@font\fontsize{#1}{#2pt}%
  \fontfamily{#3}\fontseries{#4}\fontshape{#5}%
  \selectfont}%
\fi\endgroup%
\begin{picture}(1672,2080)(261,-1835)
\put(526,-1111){\makebox(0,0)[lb]{\smash{{\SetFigFont{12}{14.4}{\rmdefault}{\mddefault}{\updefault}{\color[rgb]{0,0,0}$z$}%
}}}}
\put(526,-1786){\makebox(0,0)[lb]{\smash{{\SetFigFont{12}{14.4}{\rmdefault}{\mddefault}{\updefault}{\color[rgb]{0,0,0}$v$}%
}}}}
\put(1501,-1486){\makebox(0,0)[lb]{\smash{{\SetFigFont{12}{14.4}{\rmdefault}{\mddefault}{\updefault}{\color[rgb]{0,0,0}$w$}%
}}}}
\put(1426,-1036){\makebox(0,0)[lb]{\smash{{\SetFigFont{12}{14.4}{\rmdefault}{\mddefault}{\updefault}{\color[rgb]{0,0,0}$w_1$}%
}}}}
\put(261, 89){\makebox(0,0)[lb]{\smash{{\SetFigFont{12}{14.4}{\rmdefault}{\mddefault}{\updefault}{\color[rgb]{0,0,0}$x$}%
}}}}
\put(1571, 89){\makebox(0,0)[lb]{\smash{{\SetFigFont{12}{14.4}{\rmdefault}{\mddefault}{\updefault}{\color[rgb]{0,0,0}$y$}%
}}}}
\end{picture}%
}

Let $z$ be the last common point of
the geodesics $\overline{v x}$ and $\overline{vy}$.
We now apply Lemma~\ref{lemma:illegal:1:circuit} to the points
$\langle w_1, x, z,y \rangle$. Note that $r(w_1) < \rho_2$ (because of
$w$). Also by assumption, $r(x) \ge \rho_2 >r(w_1)$ and $r(y) \ge
\rho_2 > r(w_1)$.
Note that $h(z) = h(w_1)$, hence  $r(z) = h(z) - h(v) = h(w_1) - h(v)
> h(w_1) - h(w) = r(w_1)$. Hence we get a contradiction by
Lemma~\ref{lemma:illegal:1:circuit}. Hence we cannot have
condition $(c)$ of Lemma~\ref{lemma:family} therefore all but
$O(\eta)$ of the geodesics in $\cF'$ are unobstructed. Therefore
the number of unobstructed geodesics in $\cF'$ is at least
\begin{equation}
\label{eq:Fj:prime:many:unobstructed} (1-O(\eta)) |\cF'| \ge
(1-O(\eta)) (1-\theta_4) |E_{\lambda}| \gb^{\lambda}
\end{equation}
where we have used (\ref{eq:bounds:Fj:prime}) to get the second
estimate.

Now, let $U' \subset I_{\lambda}(v)$ be the set of $\hat{S}$-vertices (at
height $h(v)$) which are the endpoints of at least two geodesics in
$\cF'$. Since every vertex can be reached by at most $|E_{\lambda}|$
geodesics, we have by (\ref{eq:Fj:prime:many:unobstructed}),
\begin{equation}
\label{eq:Uj:prime:large:measure} |U'| \ge  (1-O(\eta)) (1-\theta_4)
\gb^{\lambda}.
\end{equation}
I.e., then $U'$ has almost full measure in $I_\lambda(v)$.

Now suppose $w \in I_\lambda(v)$ is such that two unobstructed
geodesics in $\cF'$ end at $w$. Let us denote these geodesics by
$\overline{w x}$ and $\overline{w y}$ where $x,y \in
I_\lambda'(v)$. By definition of $\cF'$, $\overline{x v}$ and
$\overline{y v}$ are unobstructed. We now apply
Lemma~\ref{lemma:illegal:1:circuit} to the points $\langle
w,x,v,y\rangle$. Note that $r(v) = 0$ (since $v$ is an
$H$-vertex), and also $r(x) \ge \rho_2$, and $r(y) \ge \rho_2$.
Thus, by Lemma~\ref{lemma:illegal:1:circuit}, we get a
contradiction unless $r(w) = 0$, i.e. $w$ is an $H$-vertex.

 \qed\medskip

If $v{\in}\cV_4$ then the conclusion is strengthened automatically
to imply that most vertices in $I_{\lambda}(v)$ are in $\cV_3$.
This is used in the following proposition.

\begin{proposition}[Zero-One Law]
\label{prop:iteration:step} Suppose $v \in \cV_4$. Suppose
$\lambda$ is such that the fraction of the edges in
$\cU(v,\lambda)$ which are in $\cE_4$ and are unobstructed for at
least length $\lambda+\rho_2$ is at least $\sigma \GG
c_9(\rho_1)/\gb^{\rho_2}$. Let $\cF = \bigcup_{w \in I_\lambda(v)}
\cU(w)$. Then at least $1 - O(\eta)$ fraction of the edges in
$\cF$ are unobstructed for length $\lambda+\rho_2$.
\end{proposition}

\bold{Proof.} As in the previous proposition, let $E_{\lambda}$ be
the set of vertices in $I'_\lambda(v)$ which are on of
$\lambda+\rho_2$ unobstructed geodesics leaving $v$. Also let $U'
\subset I_{\lambda}(v)$ and $\cF'$ be as in Proposition
\ref{prop:extension}.

Now since $v$ is not a strange vertex, the subset $U''$ of
$I_{\lambda}(v)$ consisting of $\cV_3$ vertices in $U'$ is of
almost full measure in $I_{\lambda}(v)$. Let
\begin{displaymath}
\cF'' = \bigcup_{ w \in U''} \cU(w) \cap \cE_2.
\end{displaymath}
(so $\cF''$ consists of all the $\cE_2$ edges coming out of all
the ``good'' $H$-vertices on $I_\lambda(v)$). We cut off all the
geodesics in $\cF''$ after they cross $I_\lambda'(v)$.

We want to apply Lemma \ref{lemma:family} to $\cF''$ in the box of
size $\lambda$, but there are technical difficulties here in
verifying Lemma \ref{lemma:family}.  To overcome these
difficulties, we look at a horocircle $H'$ that is $\rho_4$ units
below $H$ with the same orientation.  By a discussion similar to
Lemma \ref{lemma:comb:sameH} and following and the fact that
$\cF''$ consists of edges in $\cE_2$, we can extend every geodesic
segment in $\cF''$ by $\rho_4$ on top and bottom in all possible
ways to obtain a family $\cF'_{long}$. We will apply Lemma
\ref{lemma:family} to $\cF''_{long}$ instead. If almost all
segments in $\cF''_{long}$ are unobstructed by $H'$, it is
immediate that almost all segments in $\cF''$.  We let $U'_{long}$
be the set of $H'$ vertices within $\rho_4$ of $U'$.

 We have that $|\cF''_{long}| \ge (1-O(\eta))\gb^{2(\lambda+\rho_4)}$, so
(b) is satisfied. Also (a) is satisfied since the relevant edges
are in $\cE_2$.  If (c) does not hold, we are done, so we assume
(c) holds.  This implies that the conclusion of the lemma is true,
and we show this yields an illegal circuit. (see
Figure~\ref{fig:zero:one}).

\makefig
{Proof of Proposition~\ref{prop:iteration:step}. \\
The filled boxes denote $H'$-vertices.}{fig:zero:one}{
\begin{picture}(0,0)%
\includegraphics{circuit2.pstex}%
\end{picture}%
\setlength{\unitlength}{3947sp}%
\begingroup\makeatletter\ifx\SetFigFont\undefined%
\gdef\SetFigFont#1#2#3#4#5{%
  \reset@font\fontsize{#1}{#2pt}%
  \fontfamily{#3}\fontseries{#4}\fontshape{#5}%
  \selectfont}%
\fi\endgroup%
\begin{picture}(2323,2389)(439,-1994)
\put(1276,-1936){\makebox(0,0)[lb]{\smash{{\SetFigFont{12}{14.4}{\rmdefault}{\mddefault}{\updefault}{\color[rgb]{0,0,0}$v$}%
}}}}
\put(676,-1936){\makebox(0,0)[lb]{\smash{{\SetFigFont{12}{14.4}{\rmdefault}{\mddefault}{\updefault}{\color[rgb]{0,0,0}$u_1$}%
}}}}
\put(1876,-1936){\makebox(0,0)[lb]{\smash{{\SetFigFont{12}{14.4}{\rmdefault}{\mddefault}{\updefault}{\color[rgb]{0,0,0}$u_2$}%
}}}}
\put(526,-1261){\makebox(0,0)[lb]{\smash{{\SetFigFont{12}{14.4}{\rmdefault}{\mddefault}{\updefault}{\color[rgb]{0,0,0}$w_1$}%
}}}}
\put(2101,-1186){\makebox(0,0)[lb]{\smash{{\SetFigFont{12}{14.4}{\rmdefault}{\mddefault}{\updefault}{\color[rgb]{0,0,0}$w_2$}%
}}}}
\put(1426, 89){\makebox(0,0)[lb]{\smash{{\SetFigFont{12}{14.4}{\rmdefault}{\mddefault}{\updefault}{\color[rgb]{0,0,0}$q$}%
}}}}
\put(1426,-361){\makebox(0,0)[lb]{\smash{{\SetFigFont{12}{14.4}{\rmdefault}{\mddefault}{\updefault}{\color[rgb]{0,0,0}$q_*$}%
}}}}
\put(501,229){\makebox(0,0)[lb]{\smash{{\SetFigFont{12}{14.4}{\rmdefault}{\mddefault}{\updefault}{\color[rgb]{0,0,0}$x_1$}%
}}}}
\put(2301,239){\makebox(0,0)[lb]{\smash{{\SetFigFont{12}{14.4}{\rmdefault}{\mddefault}{\updefault}{\color[rgb]{0,0,0}$x_2$}%
}}}}
\put(1451,-1226){\makebox(0,0)[lb]{\smash{{\SetFigFont{12}{14.4}{\rmdefault}{\mddefault}{\updefault}{\color[rgb]{0,0,0}$v_j$}%
}}}}
\end{picture}%
}

By Lemma \ref{lemma:family} exists an $H'$-vertex $q$ with $h(q) <
h(v) + \lambda +2\rho_4$ and an $\hat{S}$-vertex $q_*$ with $h(q)
- \rho_2 < h(q_*) \le h(q)$ such that at least two geodesics in
$\cF''_{long}$ come together at $q_*$. Let these geodesics be
$\overline{u_1 q_*}$ and $\overline{u_2 q_*}$ where for $i=1,2$,
$u_i \in U'_{long}$. Let $w_i = \overline{u_1 q_*}{\cap}U'$ denote
the corresponding point in $U'$. Since $w_i \in U'$, there exists
$x_i \in I_{\lambda+\rho_4}'(v)$ such that $\overline{w_i x_i}$
and $\overline{x_i v}$ are both $\cE_2$ and unobstructed. Since
$v'$ denote any point on $U'_{long}$ that is $\rho_4$ units below
$v$. We now apply Lemma~\ref{lemma:illegal:1:circuit} to the
points $\langle q_*, w_1, x_1, v, x_2, w_2 \rangle$. Note that by
construction, $r(q_*) < \rho_2 \ll \rho_3$, $r(v) = \rho_4$, and
for $i=1,2$, $r(w_i) = \rho_4$, $r(x_i) \ge \rho_2$. Thus by
Lemma~\ref{lemma:illegal:1:circuit}, $q_* w_1$ and $q_* w_2$ do
not diverge at $q_*$, which is a contradiction. \qed\medskip

\begin{theorem}[Extension of Horocycles II]
\label{theorem:fullbox:align} Suppose $v \in \cV_4$ is
$x$-oriented. Let $s$ denote the height difference between $v$ and
the top of $B(L')$ and assume $s>4 \kappa^2 \beta'' R$. Then, the
density of $x$-oriented $\cV_3$ $H$-vertices along $I_s(v)$ is $1
- O(\eta)$.
%, and almost all the geodesics going up from $H$-vertices in
%$I_s(v)$ are $\cE_2$ edges which
%are unobstructed all the way to the top of $B(L')$.
\end{theorem}

\noindent{\bf Remark:} The proof of this Theorem is considerably
simpler in the case of $\DL$-graphs as boxes in $\DL$ graphs have
``no sides".  We give the proof first in this case.  The $\Sol$ case
is complicated by needing to avoid having paths ``escape off the
sides of the box."

\bold{Proof for $\DL$ graphs.}
%Fix $\sigma = O(\eta)$.
For an $x$-oriented $\cV_4$ vertex $w$, let where $f(w,\lambda)$
denotes the proportion of edges in $\cU(w)$ which are $\cE_4$ and
unobstructed for length $\lambda+\rho_2$.
%Similarly, let
%$f_j(w,\lambda)$ denote the proportion of edges in $\cP_j(w)$ which
%are $\cE_4$ and unobstructed for length $\lambda+\rho_3$.
% Clearly,
% \begin{displaymath}
% f(w,\lambda+\rho_3) = \sum_j f_j(w,\lambda+\rho_2).
% \end{displaymath}
Let
\begin{displaymath}
f^*(v,\lambda) = \sup_{w \in I_\lambda(v) \cap \cV_4}
f(w,\lambda).
\end{displaymath}
%and similarly,
%\begin{displaymath}
%f^*_j(v,\lambda) = \sup_{w \in I_\lambda(v) \cap \cV_4} f_j(w,\lambda).
%\end{displaymath}
In view of Proposition~\ref{prop:iteration:step}, for any $\lambda$
for which $f^*(v,\lambda) \ge O(\eta)$, $f^*(v,\lambda) > 1 -
O(\eta)$.

Thus, either for all $1 \le \lambda \le s$, $f^*_j(v,\lambda) \ge
1 - O(\eta)$, in which case Theorem~\ref{theorem:fullbox:align}
holds in view of Proposition~\ref{prop:extension} and
\ref{prop:iteration:step} , or else there exists minimal $\lambda$
such that $f^*(v,\lambda)
> 1-O(\eta)$, and also $f^*(v,\lambda+1) < O(\eta)$.Note that $\lambda>\Omega(\beta''R)$
by the definition of good vertices and the $\hat{S}$ and $H$-graphs.
Let $w \in I_\lambda(v) \cap \cV_4$ be such that the $\sup$ in the
definition of $f^*(v,\lambda)$ is realized at $w$.
%Note that since $\rho_2
%\ll \rho_3-1$, we have (since $\lambda$ is minimal), $f^*(v,\lambda)
%\ge 1-O(\eta)$.
Hence, by Proposition~\ref{prop:extension} (i), all but $O(\eta)$
fraction of the $\hat{S}$-vertices in $I_\lambda(w) = I_\lambda(v)$ are
$H$-vertices. By the choice of $w$ at least $1-O(\eta)$ fraction of
the geodesics in $\cU(w)$ are in $\cE_4$, unobstructed for length
$\lambda+\rho_2$, and hit an $H$-vertex (in $\cV_1$) at length
$\lambda+\rho_2+1$. Thus, in particular, the density of $H$-vertices
on $I'_{\lambda+\rho_2+1}(w)$ is at least $1 - O(\eta)$.

Let $\tilde H = I_{\lambda}(w){\cap}\cV_1$.  We consider the
family $\cE(w)$ of monotone geodesic segments ``going up" length
$L'$ from points at height $h_1$ in $\Sh(N(\phi{\inv}(\tilde
H),O(\epsilon' R)){\cap}H,\rho_1)$ and use the behavior of this
family to derive a contradiction.  We first modify $\cE(w)$ by
throwing away some bad parts of the set. This modification is
unnecessary if we are assuming that $\phi|_{U_i}$ is within
$O(\epsilon' R)$ of a $b$-standard map. We throw out any geodesic
$\gamma$ in $\cE(w)$ whose intersection with $SL^1_2(H)$ has more
than $100c_2^{\frac{1}{2}}$ of it's measure outside
$SL^1_2(H){\cap}U_*$.  By Lemma \ref{lemma:measure:lower:bilip},
this throws away at most $O(c^{\frac{1}{2}})$ of the geodesics in
$\cE(w)$. After this modifications, it follows that each geodesic
in $\cE(w)$ has $\epsilon'$-monotone image on an initial segment
of length at least $\Omega(\beta'R)$.

Note that $N(\phi{\inv}(\tilde H),O(\epsilon' R)){\cap}H$ contains a
set of large measure in $H$ and that
$(I_{\lambda}(w){\cup}I'_{\lambda+\rho_2+1}(w)){\cap}\cV_1$ is
contained in the $O(\epsilon' R)$ neighborhood of $\phi(H)$.

Since every geodesic in $\cE(w)$ diverges linearly from $H$ and
the initial segments of all $\cE(w)$ of length $\Omega({\beta'R})
> \Omega( \epsilon' R)$ have $\epsilon'$-monotone image for $\phi$,
we have that any quasi-geodesic in $\phi(\cE(w))$ diverges
linearly from $\phi(H)$ and in particular, never comes within
$\Omega(\beta' R)$ of $\phi(H)$.

Let $Q_u \subset I'_{\lambda}(w)$ be the subset of vertices $v$
such that all vertices on $I'_{\lambda+\rho_2+1}$ within
$\frac{R}{100\kappa^3}-\rho_2-1$ are not in $\cV_1$.  Since
$\ell(\cV_1^c{\cap}I'_{\lambda+\rho_2+1})<O(\eta)\ell(
I'_{\lambda+\rho_2+1})$, we have $\ell(Q_u) \ll
O(\eta)\ell(I'_{\lambda}(w))=O(\eta)\gb^{\lambda}$. Any
quasi-geodesic in $\phi(\cE(w))$ crossing $I'_{\lambda}(w)$ does
so on $Q_u$.

Similarly let $Q_d = I_{\lambda}(w){\cap}\cV_1^c$ and note that
$\ell(Q_d) \ll O(\eta)\ell(I_{\lambda}(w))=O(\eta)\gb^{\lambda}$
and that any quasi-geodesic in $\phi(\cE(w))$ crossing
$I_{\lambda}(w)$ must cross it on $Q_d$.

Now as all quasi-geodesics in $\phi(\cE(w))$ diverge linearly from
$\phi(H)$, they must all eventually leave the box of size $\lambda$
bounded by $I_{\lambda}(w)$ and $I'_{\lambda}(w)$.  This implies
that every quasi-geodesic in $\phi(\cE(w))$ eventually crosses
$Q_u{\cup}Q_d$ or that every geodesic in $\cE(w)$ eventually crosses
$\phi{\inv}(Q_u{\cup}Q_d)$.  This is impossible by Lemma
\ref{lemma:length:area}, since
\begin{equation}
\label{eq:Qu:Qd:small}
\ell(\phi^{-1}({Q_u} \cup {Q_d})) \le O(\eta) \ell(H),
\end{equation}
and $c(\rho_1) O(\eta)  \ll 1$. \qed\medskip

Before reading the proof for $\Sol$, the reader should be sure to
read \S\ref{subsection:tangle}.

\bold{Proof for $\Sol$.}
We need to modify the proof given above in two ways in order to
avoid ``escape off the sides" of the box of size $\lambda$.  As
this is a modification of the previous proof, we only sketch the
necessary changes.

% Given an
%$x$-oriented (resp. $y$ oriented) horocycle $H$ in either the domain
%or range, we let $A(H,l)=\Sh(H,\rho_1)\cap h{\inv}(h(H)+l)$ (resp.
%$\Sh(H,\rho_1)\cap h{\inv}(h(H)-l)$.

We choose $w$ as in the proof for $DL$ graphs.  We remark that it
is easy to see that $w$ can be chosen away from the edge of
$B[\lambda]$. This can be deduced from Proposition
\ref{prop:iteration:step}.  We will assume that we have chosen
such a $w$.  It is also possible to work with $w$ near the edge of
the box, but that one use a more complicated definition of points
deep in the shadow of horocycles.

As above we consider $\tilde H = I_{\lambda}(w){\cap}\cV_1$. We
consider the family $\cE(w)$ of monotone geodesics ``going up"
length $L'$ from points at height
$$h_1 \in \Sh(N(\phi{\inv}(\tilde H),O(\epsilon' R)){\cap}H,\rho_1)$$
and
use the behavior of this family to derive a contradiction. We
first modify $\cE(w)$ exactly as before. We now further modify
$\cE(w)$ to only include those geodesics whose images at the end
of the initial segment are $\beta'R$-deep in $B[\lambda]$. By this
we mean that they are $\beta'R$ deep in the shadows of the top and
bottom of $B(R)$. This subset still contains a large proportion of
the original elements of $\cE(w)$. Let $Q=Q_u \cup Q_d$.  Then as
before, we see that paths in $\cE(w)$ can only come near the top
and bottom of $B[\lambda]$ in $N(Q^c, \frac{\beta'R}{2 \kappa})$.

We now apply the results of \S\ref{subsection:tangle} with
$\rho=\rho_1, D_1=\epsilon'R, D_2=\frac{\beta'R}{2\kappa}$ and
$D_3=\lambda$. By Lemma~\ref{lemma:shadows:tangle} if a path
$\gamma \in \phi(\cE(w))$ leaves the box, it must tangle with the
union of the top and the bottom of the box.  Since $\gamma \in
N(Q^c, \frac{\beta'R}{2 \kappa}))$, Lemma
\ref{lemma:generalized:trapping} implies
$$\ell(Q_u \cup Q_d)
=\ell(Q) \ge \omega \| \cE(w)\|,$$ \noindent where $\omega$
depends only on $\kappa$ and $C$. But we have $\|\cE(w)\| \ge
\omega' \ell(H)$. where $\omega'$ depends only on $\kappa$ and
$C$. This is a conradiction to (\ref{eq:Qu:Qd:small}), if $\eta$
is sufficiently small. As before $\eta$ can be made arbitrarily
small by taking $\epsilon'$ and $\delta_0$ sufficiently small.
\qed\medskip

\subsection{Completing the proof of Theorem
  \ref{theorem:main:step:two}}
\label{sec:subsection:end:proof:theorem2.1}

\bold{Proof.}  By Theorem~\ref{theorem:lamplighter} $\phi{\inv}$
of any very favorable horocycle  in $B(L')$ is within $O(\epsilon'
R)$ error of a horocycle.  Given $\hat \theta>0$,  Lemma
\ref{lemma:most:favorable} implies that, by choosing $\beta''$ and
$\delta_0$ small enough, that $1-\hat \theta$ of the measure of
$B(L')$ consists of points in the image of both a very favorable
$x$-horocycle and a very favorable $y$-horocycle.  By an argument
from the proof of \cite[Lemma 4.11]{EFW1}, this implies that
$\phi{\inv}$ respects level sets of height to within $O(\epsilon'
R)$ error.

\relax From this, it is not hard to show that $\phi{\inv}$ of most
vertical geodesics are weakly monotone.  This is very similar to the
proof of \cite[Lemma 6.5]{EFW1}.  There are some additional
difficulties due to the fact that we only control the map on most of
the measure, but these can be handled in a manner similar to the
proofs of \cite[Lemma 5.10 and Corollary 5.12]{EFW1}.

Once we know $\phi{\inv}$ of most vertical geodesics are weakly
monotone, the conclusion of the theorem follows from as in the
proof of \cite[Theorem 5.1]{EFW1} at the end of \cite[Section
5.4]{EFW1}. \qed

\noindent  Department of Mathematics, University of Chicago. \\
5734 S. University Avenue, Chicago, Illinois 60637.\

\medskip
\noindent Department of Mathematics, Indiana University. \\
Rawles Hall, Bloomington, IN, 47405.\

\medskip
\noindent Department of Mathematics, Statistics, \& Computer
Science, University of Illinois at Chicago. 322 Science
Engineering Offices (M/C 249), 851 S. Morgan Street Chicago, IL
60607-7045.\


\medskip


\begin{thebibliography}{Palais}

%\bibitem[BLPS]{BLPS}
%Benjamini, I.; Lyons, R.; Peres, Y.; Schramm, O. Group-invariant
%percolation on graphs.  {\it Geom. Funct. Anal.}  9  (1999),  no.
%1, 29--66.

%\bibitem[Co]{Co} Cooper, Daryl, Appendix to \cite{FM1}.

%\bibitem[Ch]{Ch}
%Chow, Richard Groups quasi-isometric to complex hyperbolic space.
%{\it Trans. Amer. Math. Soc.} 348 (1996), no. 5, 1757--1769.

\bibitem[dlH]{dlH}
de la Harpe, Pierre. {\it Topics in geometric group theory.}
Chicago Lectures in Mathematics. University of Chicago Press,
Chicago, IL, 2000.

%\bibitem[DL]{DL}
%Diestel, Reinhard; Leader, Imre A. conjecture concerning a limit
%of non-Cayley graphs.  {\it J. Algebraic Combin.}  14  (2001), no.
%1, 17--25.

%\bibitem[Dy]{Dy}
%Dymarz, Tullia. Ph.d. Thesis, University of Chicago, 2006.

%\bibitem[D]{D}
%Dyubina, Anna. Instability of the virtual solvability and the
%property of being virtually torsion-free for quasi-isometric
%groups.  {\em Internat. Math. Res. Notices}  2000,  no. 21,
%1097--1101.

%\bibitem[ET]{ET}
%Elek, Gabor; Tardos, Gabor. On roughly transitive amenable graphs
%and harmonic Dirichlet functions. {\it Proc. Amer. Math. Soc.} 128
%(2000), no. 8, 2479--2485.


%\bibitem[E]{E} Eskin, Alex. Quasi-isometric rigidity of
%nonuniform lattices in higher rank symmetric spaces. {\it J. Amer.
%Math. Soc.}  11 (1998), no. 2, 321--361.

%\bibitem[EF]{EF}
%Eskin, Alex; Farb, Benson. Quasi-flats and rigidity in higher rank
%symmetric spaces. {\it J. Amer. Math. Soc. } 10  (1997),  no. 3,
%653--692.

\bibitem[EFW1]{EFW0}
Eskin, Alex; Fisher, David; Whyte, Kevin. Quasi-isometries and
rigidity of solvable groups, {\it Pure Appl. Math. Q.} 3 (2007), no. 4, part 1, 927–947.


\bibitem[EFW2]{EFW1}
Eskin, Alex; Fisher, David; Whyte, Kevin. Coarse Differentiation
of Quasi-isometries I; spaces not quasi-isometric to Cayley
graphs, to appear {\it Annals of Math.}


%\bibitem[FS]{FS}
%Farb, Benson; Schwartz, Richard. The large-scale geometry of
%Hilbert modular groups. {\it J. Differential Geom.}  44  (1996),
%no. 3, 435--478.

%\bibitem[F]{F} Farb, Benson. The quasi-isometry classification of lattices
%in semisimple Lie groups.  {\it Math. Res. Lett.}  4  (1997),  no.
%5, 705--717.

%\bibitem[FM1]{FM1} Farb, Benson; Mosher, Lee. A rigidity theorem for the
%solvable Baumslag-Solitar groups. With an appendix by Daryl
%Cooper. {\it Invent. Math.}  131  (1998),  no. 2, 419--451.

%\bibitem[FM2]{FM2}
%Farb, Benson; Mosher, Lee. Quasi-isometric rigidity for the
%solvable Baumslag-Solitar groups. II. {\it Invent. Math.}  137
%(1999),  no. 3, 613--649.

%\bibitem[FM3]{FM3}
%Farb, Benson; Mosher, Lee. On the asymptotic geometry of
%abelian-by-cyclic groups.  {\it Acta Math.}  184  (2000),  no. 2,
%145--202.

\bibitem[FM]{FM4}
Farb, Benson; Mosher, Lee. Problems on the geometry of finitely
generated solvable groups. {\it  Crystallographic groups and their
generalizations (Kortrijk, 1999)} ,  121--134, Contemp. Math.,
262, Amer. Math. Soc., Providence, RI, 2000.

\bibitem[GhdlH]{GhdlH}
Sur les groupes hyperboliques d'apr\'{e}s Mikhael Gromov. (French)
[Hyperbolic groups in the theory of Mikhael Gromov] Papers from
the Swiss Seminar on Hyperbolic Groups held in Bern, 1988. Edited
by \'{E}. Ghys and P. de la Harpe. Progress in Mathematics, 83.
Birkhäuser Boston, Inc., Boston, MA, 1990.

%\bibitem[Gr1]{Gr1}
%Gromov, Mikhael. Groups of polynomial growth and expanding maps.
%{\it Inst. Hautes \'{E}tudes Sci. Publ. Math.} No. 53 (1981),
%53--73.

%\bibitem[Gr2]{Gr2}
%Gromov, Mikhael. Infinite groups as geometric objects. Proceedings
%of the International Congress of Mathematicians, Vol. 1, 2
%(Warsaw, 1983),  385--392, PWN, Warsaw, 1984.

%\bibitem[Gr3]{Gr3}Gromov, M. Asymptotic invariants of infinite groups.
%{\it Geometric group theory, Vol. 2 (Sussex, 1991)},  1--295,
%London Math. Soc. Lecture Note Ser., 182, Cambridge Univ. Press,
%Cambridge, 1993.

%\bibitem[He]{He}
%Heinonen, Juha {\it Lectures on analysis on metric spaces.}
%Universitext. Springer-Verlag, New York, 2001.

%\bibitem[H]{H}
%Hinkkanen, A. Uniformly quasisymmetric groups. {\it Proc. London
%Math. Soc.} (3) 51 (1985), no. 2, 318--338.

%\bibitem[KL]{KL}
%Kleiner, Bruce; Leeb, Bernhard. Rigidity of quasi-isometries for
%symmetric spaces and Euclidean buildings. {\it Inst. Hautes
%\'{E}tudes Sci. Publ. Math.} No. 86, (1997), 115--197 (1998).

%\bibitem[K]{K} Kleiner, Bruce. Personal communication.

%\bibitem[KR]{KR} Kor�nyi, A.; Reimann, H. M. Foundations for the theory of
%quasiconformal mappings on the Heisenberg group.  Adv. Math.  111
%(1995),  no. 1, 1--87.

\bibitem[MN]{MN}
Letter from R.Moeller to W.Woess, 2001.

%\bibitem[MSW]{MSW}
%Mosher, Lee; Sageev, Michah; Whyte, Kevin. Quasi-actions on trees.
%I. Bounded valence.  {\it Ann. of Math.} (2)  158  (2003),  no. 1,
%115--164.

%\bibitem[Mo1]{Mo1}
%Mostow, G. D. Quasi-conformal mappings in $n$-space and the
%rigidity of hyperbolic space forms.  Inst. Hautes \'{E}tudes Sci.
%Publ. Math. No. 34 1968 53--104.

%\bibitem[Mo2]{Mo0}
%Mostow, G. D. Representative functions on discrete groups and
%solvable arithmetic subgroups. {\it Amer. J. Math.} 92 1970 1--32.

%\bibitem[Mo2]{Mo2}
%Mostow, G. D. Strong rigidity of locally symmetric spaces. Annals
%of Mathematics Studies, No. 78. Princeton University Press,
%Princeton, N.J.; University of Tokyo Press, Tokyo, 1973.


%\bibitem[P1]{P}
%Pansu, Pierre. Metriques de Carnot-Carath\'{e}odory et
%quasiisom\'{e}tries des espaces sym\'{e}triques de rang un.
%(French) [Carnot-Caratheodory metrics and quasi-isometries of
%rank-one symmetric spaces] {\it Ann. of Math.} (2)  129  (1989),
%no. 1, 1--60.

%\bibitem[P2]{P2}
%Pansu, Pierre. Dimension conforme et sph\`{e}re l'infini des
%vari\'{e}t\'{e}s \`{a} courbure n\'{e}gative. (French) [Conformal
%dimension and sphere at infinity of manifolds of negative
%curvature]  Ann. Acad. Sci. Fenn. Ser. A I Math.  14 (1989),  no.
%2, 177--212.

\bibitem[PPS]{PPS}
Yuval Peres, Gabor Pete, Ariel Scolnicov. Critical percolation on
certain non-unimodular graphs, {\it New York J. Math.} 12 (2006), 1–18 (electronic).

%\bibitem[Sa]{Sa}
%Sauer, Roman. Homological Invariants and Quasi-Isometry. To appear
%{\it GAFA}.

%\bibitem[S1]{S1}
%Schwartz, Richard Evan. The quasi-isometry classification of rank
%one lattices.  {\it Inst. Hautes \'{E}tudes Sci. Publ. Math.} No.
%82 (1995), 133--168 (1996).

%\bibitem[S2]{S2}
%Schwartz, Richard Evan. Quasi-isometric rigidity and Diophantine
%approximation.  {\it Acta Math.}  177  (1996),  no. 1, 75--112.

%\bibitem[Sh]{Sh}
%Shalom, Yehuda. Harmonic analysis, cohomology, and the large-scale
%geometry of amenable groups. {\it Acta Math.}  192  (2004),  no.
%2, 119--185.

%\bibitem[SW]{SW}
%Soardi, Paolo M.; Woess, Wolfgang. Amenability, unimodularity, and
%the spectral radius of random walks on infinite graphs. {\it Math.
%Z.} 205 (1990), no. 3, 471--486.

%\bibitem[Tu]{Tu}
%Tukia, Pekka. On quasiconformal groups.  {\it J. Analyse Math.} 46
%(1986), 318--346.

%\bibitem[Wo1]{Wo}
%Woess, Wolfgang. Topological groups and infinite graphs. {\it
%Directions in infinite graph theory and combinatorics (Cambridge,
%1989). Discrete Math.} 95 (1991), no. 1-3, 373--384.

\bibitem[Wo]{Wo2}
Woess, Wolfgang.  Lamplighters, Diestel-Leader graphs, random
walks, and harmonic functions, {\it Combinatorics, Probability \&
Computing} 14 (2005) 415-433.

\bibitem[W]{W}
Wortman, Kevin. A finitely presented solvable group with small
quasi-isometry group, {\it Michigan Math. J.} 55 (2007), no. 1, 3–24.

\end{thebibliography}
\end{document}